# Rough Fuzzy Quadratic Minimum Spanning Tree Problem

Saibal Majumder[1], Samarjit Kar[2], Tandra Pal[3]

[1,3] Department of Computer Science, National Institute of Technology, Durgapur, West Bengal, India- 713209
[2] Department of Mathematics, National Institute of Technology, Durgapur, West Bengal, India- 713209
e-mail: saibaltufts@gmail.com; kar_s_k@yahoo.com

**Abstract:** A quadratic minimum spanning tree (QMST) problem is to determine a minimum spanning tree of a connected graph having edges which are associated with linear and quadratic weights. The linear weights are the edge costs which are associated with every edge whereas the quadratic weights are the interaction costs between a pair of edges of the graph. In this paper, a bi-objective $(\alpha, \beta)$ rough fuzzy quadratic minimum spanning tree problem (b-$(\alpha, \beta)$RFQMSTP) has been considered for a connected graph whose linear and quadratic weights are expressed as rough fuzzy variable. $\alpha$ and $\beta$ determine the predefined confidence levels of the model. The b-$(\alpha, \beta)$ RFQMSTP is transformed into crisp form using chance-constrained programming technique. The crisp equivalent of b-$(\alpha, \beta)$ RFQMSTP model is solved using epsilon-constraint method and two multi-objective evolutionary algorithms: nondominated sorting genetic algorithm II (NSGA-II) and multi-objective cross generational elitist selection, heterogeneous recombination and cataclysmic mutation (MOCHC) algorithm. A numerical example is provided to illustrate the proposed model when solved with different techniques. A sensitivity analysis of the example is performed at different confidence levels of $\alpha$ and $\beta$. Finally a performance analysis of NSGA-II and MOCHC are performed on 5 randomly generated bi-objective $(\alpha, \beta)$ rough fuzzy quadratic minimum spanning tree instances. The result of our proposed model when solved by epsilon-constraint method is obtained by using the standard optimization software, LINGO. Whereas the results of evolutionary algorithms for our model are obtained by using an open source optimization package, jMETAL [22].

**1. Introduction:** Minimum spanning tree (MST) problem is one of the most fundamental and important problems of graph theory which has diverse application in different engineering and scientific domains. A MST problem is to find a minimum cost spanning tree of a weighted connected network (WCN). In the literature the problem has been well studied by Prime [33], Kruskal [19], Graham and Hell [35], Gabow et al. [17], Bondy and Murty [18] and Bazlamaçi and Hindi [15]. QMST problem is a NP-hard problem which is a variant of MST introduced by Assad and Xu [1]. The problem is to find a minimum spanning tree of a WCN which includes two types of weights; linear and quadratic. The linear weights are associated with every edge. And the quadratic weights, also known as intercosts, are due to the effect of interaction between any pair of edges. QMST problem has applications in different real networks. As an example, in telecommunication network [34] the linear costs represent the cost incurred to establish a radio or cable connection while the quadratic weights between a pair of connected edges represent the cost of conversion device between different types of cables. In a utility network (oil or water transmission network) or in transportation network the linear weights represent the cost to install different pipes or to or construct roads between two junctions. Moreover, the quadratic weights represent the interaction/interfacing costs for installing valves or bending losses in T-joints [34] in utility network, and turn penalties in transportation network [1]. These interfacing costs mostly exist for adjacent edges in these cases. However, sometimes these interfacing cost crops up for any pair of connected edges. But such topological design has little resemblance with physical layout, as in the case of fibre-optic networks [34].

In the literature there are some studies related to QMST. Assad and Xu [1] provided a branch and bound based exact method for the problem. The result of the QMST obtained by [1] was improved by a genetic algorithm by Zhou and Gen [16]. Sundar and Singh [36] provided an artificial bee colony algorithm for the problem. Cordone and Passeri [34] developed a tabu search and a variable neighbourhood search for the problem. Maia et al. [37] proposed a Pareto local search algorithm for a bi-objective adjacent only QMST, where the quadratic weights are associated with a pair of adjacent edges only. All these related studies have been done when the related problem is defined under precise and crisp ambiance.

In many real world scenarios it becomes difficult to define the exact linear and quadratic weights for the QMST due to the uncertainties associated with the weights. The uncertain nature of these weights may be due to insufficient information or vagueness, error in measurements, etc. Under such uncertain environment there are two existing works in the literature. In 2005, Gao and Lu [20] defined a fuzzy QMST and solved the problem using a genetic algorithm. The fuzzy weights associated with the problem can only process the uncertainty represented in terms of impreciseness. Recently in 2014, Zhou et al. [26] defined a chance-constrained model for a QMST using Liu's uncertainty theory [7]. The uncertainty theory [7] is based on modelling human belief degree estimated by human/experts, i.e., the strength with which we believe the event will occur. However, there are many real life problems where situation demands a mixture of uncertainty in terms of impreciseness and vagueness. In such cases rough fuzzy theory [5, 6] can be more effective rather than fuzzy and uncertainty theory which are separately used in the existing studies [20, 26] to model QMST problem in uncertain environment. In practice, different decision makers (DMs) or experts have different opinions about the linear and quadratic weights for a WCN. These weights may be expressed as intervals. Moreover these interval weights may fluctuate from time to time, e.g., for a road network, the linear and quadratic interval costs may fluctuate over a time span due to fluctuation in fuel price, toll tax or turn penalties. Specifically, the expenses related with turn penalties of a road is directly proportional with fuel price and the traffic congestion in the road at a particular time; hence the uncertain information about these weights may be a mixture of both impreciseness and vagueness. Different types of such uncertainties may well be expressed as rough fuzzy variables [5, 6]. In the literature there exists no study on QMST under such mixed uncertain environment which motivates us to study rough fuzzy bi-objective QMST. We have defined both the linear and quadratic weights as rough fuzzy variables [5, 6]. The bi-objective QMST has been modelled by rough fuzzy chance-constrained programming (RFCCP) technique. The linear weights and quadratic weights are used as the two objectives of the model. The model is solved by epsilon-constraint method [39] and two multi-objective evolutionary algorithms (MOEAs): NSGA-II [27] and MOCHC [23]. Finally a sensitivity analysis at different confidence level is performed for the RFCCP model of QMST.

The rest of the paper is organized as follows. Section 2 provides some basic concepts related to uncertainty like fuzzy variables, rough variables, rough fuzzy variables, etc. Section 3 introduces the concepts of multi-objective optimization. Section 4 presents some results of RFCCP model. The proposed chance-constrained model is designed in Section 5. Section 6 discusses epsilon-constraint method, NSGA-II [27] and MOCHC [23] algorithms. The RFCCP model of b-$(\alpha, \beta)$RFQMSTP is numerically illustrated by three methods and the result are discussed in Section 7. Finally Section 8 culminates the paper.

## 2. Preliminaries

In this section we revise the concepts of fuzzy variables, rough variables, rough fuzzy variables and some related properties.

**Definition 1 [5]:** A fuzzy variable $\tilde{\xi}$ is defined as a function from a possibility space $(\Theta, \mathbb{P}(\Theta), Pos)$ to a set of real numbers $\Re$, where $\Theta$ represents a nonempty set of sample space, $\mathbb{P}(\Theta)$ is the power set of $\Theta$ and $Pos$ is the possibility measure of a fuzzy event $\{\tilde{\xi} \in A\}, A \in \Re$, which is defined as

$$Pos\{\tilde{\xi} \in A\} = sup_{x \in A}\tilde{\xi}(x) \tag{1}$$

A fuzzy variable $\tilde{\xi}$ is said to be a triangular fuzzy variable (TFV) if it is represented by $\tilde{\xi} = (\xi_1, \xi_2, \xi_3)$, where, $\xi_1 < \xi_2 < \xi_3$ and $\tilde{\xi}$ follows a distribution for its membership function $\mu_{\tilde{\xi}}(x \in \Re)$ as given in (2)

$$\mu_{\tilde{\xi}}(x \in \Re) = \begin{cases} 0 & ; x < \xi_1 \\ \frac{x-\xi_1}{\xi_2-\xi_1} & ; \xi_1 \leq x \leq \xi_2 \\ \frac{\xi_3-x}{\xi_3-\xi_2} & ; \xi_2 \leq x \leq \xi_3 \end{cases} \tag{2}$$

**Definition 2 [6]:** For a possibility space $(\Theta, \mathbb{P}(\Theta), Pos)$ the possibility of occurrence of a fuzzy event, $\mathcal{F} \in \mathbb{P}(\Theta)$ is represented as $Pos\{\mathcal{F}\}$ while the necessity of occurrence of the fuzzy event $\mathcal{F}$ is defined as follows in (3).

$$Nec\{\mathcal{F}\} = 1 - Pos\{\mathcal{F}^c\} \tag{3}$$

**Definition 3 [6]:** For a possibility space $(\Theta, \mathbb{P}(\Theta), Pos)$, $\mathcal{F}$ be a fuzzy event $\mathcal{F} \in \mathbb{P}(\Theta)$ then the credibility measure of $\mathcal{F}$ is expressed as

$$Cr\{\mathcal{F}\} = \frac{1}{2}(Pos\{\mathcal{F}\} + Nec\{\mathcal{F}\}) \tag{4}$$

It can be observed that a fuzzy event may be rejected even though its possibility of occurrence is 1 and may hold even if its possibility is 0. However a fuzzy event must be accepted if its credibility is 1 and rejected if its credibility is 0.

**Definition 4 [6]:** The credibility distribution of a TFV $\tilde{\xi} = (\xi_1, \xi_2, \xi_3)$ can be defined as

$$Cr\{\tilde{\xi} \leq x\} = \begin{cases} 0 & ; x < \xi_1 \\ \frac{x-\xi_1}{2(\xi_2-\xi_1)} & ; \xi_1 \leq x \leq \xi_2 \\ \frac{x+\xi_3-2\xi_2}{2(\xi_3-\xi_2)} & ; \xi_2 \leq x \leq \xi_3 \\ 1 & ; x > \xi_3 \end{cases} \tag{5}$$

and

$$Cr\{\tilde{\xi} \geq x\} = \begin{cases} 1 & ; x < \xi_1 \\ \frac{2\xi_2-\xi_1-x}{2(\xi_2-\xi_1)} & ; \xi_1 \leq x \leq \xi_2 \\ \frac{\xi_3-x}{2(\xi_3-\xi_2)} & ; \xi_2 \leq x \leq \xi_3 \\ 0 & ; x > \xi_3 \end{cases} \tag{6}$$

**Lemma 2.1 [20]:** Consider $\tilde{\xi} = (\xi_1, \xi_2, \xi_3, \xi_4)$ be a trapezoidal fuzzy variable (TrFV) and $\alpha$ be the predetermined confidence level then,

i. While $\alpha \in (0, 0.5]$, $Cr\{\xi \leq x^*\} \geq \alpha$ if and only if $x^* \geq (1-2\alpha)\xi_1 + 2\alpha\xi_2$
ii. While $\alpha \in (0.5, 1]$, $Cr\{\xi \leq x^*\} \geq \alpha$ if and only if $x^* \geq (2-2\alpha)\xi_3 + (2\alpha-1)\xi_4$

**Remark 2.1:** If $\tilde{\xi} = (\xi_1, \xi_2, \xi_3)$ is a triangular fuzzy variable (TFV), which is a special case of a trapezoidal fuzzy number, then we can have,

i. While $\alpha \in (0, 0.5]$, $Cr\{\xi \leq x^*\} \geq \alpha$ if and only if $x^* \geq (1 - 2\alpha)\xi_1 + 2\alpha\xi_2$

ii. While $\alpha \in (0.5, 1]$, $Cr\{\xi \leq x^*\} \geq \alpha$ if and only if $x^* \geq (2 - 2\alpha)\xi_2 + (2\alpha - 1)\xi_3$

**Definition 5** [5, 6]: Let $\Lambda$ be a nonempty subset, $\mathcal{A}$ be a $\sigma$-algebra of subsets of $\Lambda$, $\Delta$ be an element in $\mathcal{A}$ and $\pi$ be the nonnegative real-valued additive set function, then $(\Lambda, \mathcal{A}, \Delta, \pi)$ is said to be a rough space.

**Definition 6** [5, 6]: A rough variable $\xi$ on a rough space $(\Lambda, \mathcal{A}, \Delta, \pi)$ is defined as a measurable function from $\Lambda$ to a set of real numbers $\Re$ such that for any Borel set $\mathcal{S}$ of $\Re$ we have

$$\{\lambda \in \Lambda | \xi(\lambda) \in \mathcal{S}\} \in \mathcal{A} \tag{7}$$

The lower approximation $\underline{\xi}$ and upper approximation $\overline{\xi}$ of a rough variable $\xi$ are respectively defined as

$$\underline{\xi} = \{\xi(\lambda) | \lambda \in \Delta\} \text{ and } \overline{\xi} = \{\xi(\lambda) | \lambda \in \Lambda\} \tag{8}$$

From (4), it is conspicuous that $\underline{\xi} \subset \overline{\xi}$ since $\Delta \subset \Lambda$.

**Definition 7** [5, 6]: Let $\xi$ be a $n$ dimensional rough vector on a rough space $(\Lambda, \mathcal{A}, \Delta, \pi)$ and $f_j: \Re^n \to \Re$ be the continuous functions, $j = 1, 2, \ldots, n$. Then the lower trust and upper trust of a rough event characterized by $f_j(\xi) \leq 0; j = 1, 2, \ldots, n$ are respectively defined as

$$T\underline{r}\{f_j(\xi) \leq 0; j = 1, 2, \ldots, n\} = \frac{\pi\{\lambda \in \Delta | f_j(\xi) \leq 0; j = 1, 2, \ldots, n\}}{\pi\{\Delta\}} \tag{9}$$

and

$$T\overline{r}\{f_j(\xi) \leq 0; j = 1, 2, \ldots, n\} = \frac{\pi\{\lambda \in \Lambda | f_j(\xi) \leq 0; j = 1, 2, \ldots, n\}}{\pi\{\Lambda\}} \tag{10}$$

If $\pi\{\Delta\} = 0$, then $T\underline{r}\{f_j(\xi) \leq 0; j = 1, 2, \ldots, n\} \equiv T\overline{r}\{f_j(\xi) \leq 0; j = 1, 2, \ldots, n\}$

The trust of a rough event is defined as a mean of lower and upper trusts of the same event, i.e.,

$$Tr\{f_j(\xi) \leq 0; j = 1, 2, \ldots, n\} = \frac{1}{2}\left(T\underline{r}\{f_j(\xi) \leq 0; j = 1, 2, \ldots, n\} + T\overline{r}\{f_j(\xi) \leq 0; j = 1, 2, \ldots, n\}\right) \tag{11}$$

**Example:** Let $\xi$ be a rough variable defined on a rough space $(\Lambda, \mathcal{A}, \Delta, \pi)$ such that $\xi = [\xi_1, \xi_2][\xi_3, \xi_4]$, $\xi_3 \leq \xi_1 \leq \xi_2 \leq \xi_4$, where $[\xi_1, \xi_2]$ and $[\xi_3, \xi_4]$ are respectively the lower approximation and upper approximation of $\xi$, i.e., all the elements of the interval $[\xi_1, \xi_2]$ are certainly the members of $\xi$ and the elements of the interval $[\xi_3, \xi_4]$ are the possible members of $\xi$. Thus, $\Delta = \{t | \xi_1 \leq t \leq \xi_2\}$ and $\Lambda = \{t | \xi_3 \leq t \leq \xi_4\}$, $\xi(t) = t \ \forall \ t \in \Lambda$. $\mathcal{A}$ is the $\sigma$-algebra on $\Lambda$ and $\pi$, the Lebesgue measure.

As a practical example, consider a case where one needs to travel from a source $s$ to destination $d$. The cost incurred during the journey from $s$ to $d$ depend on many factors like, cost of fuel, toll tax, conveyance cost, etc. These cost parameters are dynamic in nature since they may vary with time. Hence it is very difficult to estimate these costs in exact figures. Let four experts predict the unit cost, $c_{sd}$ of travelling from destination $s$ to $d$, as intervals: $[25, 40.8], [30.2, 42], [20.5, 50.5]$ and $[27, 45]$ respectively. From such opinion it is clear that every elements in the interval $\underline{c_{sd}} = [30.2 \ 40.8]$ are also contained in four different intervals defined by experts predictions. Hence the elements of $\underline{c_{sd}}$ are certainly regarded as the greatest definable set contained in $c_{sd}$, and is considered as the lower approximation of $c_{sd}$. Again the interval $\overline{c_{sd}} = [20.5, 50.5]$ becomes the least definable set containing $c_{sd}$ since the elements of $\overline{c_{sd}}$ partly belong to the intervals defined by experts predictions and hence is regarded as upper approximation of $c_{sd}$. Then, $c_{sd}$ is expressed as a rough variable, $[30.2, 40.8][20.5, 50.5]$. Similar examples for rough variables can also be observed in the study presented in [41].

**Definition 8** [4, 5]: The trust distribution of a rough variable $\xi = [\xi_1, \xi_2][\xi_3, \xi_4]$ is defined by (12) and (13).

$$Tr\{\xi \leq x\} = \begin{cases} 0 & ; x \leq \xi_3 \\ 0.5\left(\frac{x-\xi_3}{\xi_4-\xi_3}\right) & ; \xi_3 \leq x \leq \xi_1 \\ 0.5\left(\frac{x-\xi_1}{\xi_2-\xi_1} + \frac{x-\xi_3}{\xi_4-\xi_3}\right) & ; \xi_1 \leq x \leq \xi_2 \\ 0.5\left(\frac{x-\xi_3}{\xi_4-\xi_3} + 1\right) & ; \xi_2 \leq x \leq \xi_4 \\ 1 & ; x \geq \xi_4 \end{cases} \quad (12)$$

and

$$Tr\{\xi \geq x\} = \begin{cases} 0 & ; x \geq \xi_4 \\ 0.5\left(\frac{\xi_4-x}{\xi_4-\xi_3}\right) & ; \xi_2 \leq x \leq \xi_4 \\ 0.5\left(\frac{\xi_4-x}{\xi_4-\xi_3} + \frac{\xi_2-x}{\xi_2-\xi_1}\right) & ; \xi_1 \leq x \leq \xi_2 \\ 0.5\left(1 + \frac{\xi_4-x}{\xi_4-\xi_3}\right) & ; \xi_3 \leq x \leq \xi_1 \\ 1 & ; x \leq \xi_3 \end{cases} \quad (13)$$

**Definition 9** [5, 6]: A rough fuzzy variable $\hat{\xi}$ is a function from a possibility space $(\Theta, \mathbb{P}(\Theta), Pos)$ to a set of rough variables such that $Tr\{\hat{\xi}(\theta \in \Theta) \in S\}$ is a fuzzy variable.

**Definition 10** [5, 6]: Let $\hat{\xi}$ be a $n$ dimensional rough fuzzy vector then $\hat{\xi}$ is a function from possibility space $(\Theta, \mathbb{P}(\Theta), Pos)$ to a set of $n$ dimensional rough vectors such that $Tr\{\hat{\xi}(\theta \in \Theta) \in S\}$ is a measurable map of $\theta$ for any Borel set $S$ of $\Re^n$.

**Definition 11** [5, 6]: Let $u: \Re^n \to \Re$ be a function and $\hat{\xi}_1, \hat{\xi}_2, \ldots, \hat{\xi}_n$ be the rough fuzzy variables on possibility space $(\Theta, \mathbb{P}(\Theta), Pos)$, then $\hat{\xi} = u(\hat{\xi}_1, \hat{\xi}_2, \ldots, \hat{\xi}_n)$ is a rough fuzzy variable and determined by, $\hat{\xi}(\theta) = u(\hat{\xi}_1(\theta), \hat{\xi}_2(\theta), \ldots, \hat{\xi}_n(\theta)) \,\forall\, \theta \in \Theta$ \hfill (15)

**Definition 12** [5, 6]: Let $u: \Re^n \to \Re$ is a continuous function and $\xi_i$ be the rough fuzzy variable defined on $(\Theta_i, \mathbb{P}(\Theta_i), Pos_i) \,\forall i \in \{1, 2, \ldots, n\}$. Then $\hat{\xi} = u(\hat{\xi}_1, \hat{\xi}_2, \ldots, \hat{\xi}_n)$ is a rough fuzzy variable defined on the product possibility space $(\Theta, \mathbb{P}(\Theta), Pos)$ such that,
$$\hat{\xi}(\theta_1, \theta_2, \ldots, \theta_n) = u(\hat{\xi}_1(\theta_1), \hat{\xi}_2(\theta_2), \ldots, \hat{\xi}_n(\theta_n)) \,\forall\, \theta_i \in \Theta. \quad (16)$$

**Definition 13** [5, 6]: Let $u_j: \Re^n \to \Re, j = 1, 2, \ldots, m$ are continuous functions and $\hat{\xi} = (\hat{\xi}_1, \hat{\xi}_2, \ldots, \hat{\xi}_n)$ be the rough fuzzy vector defined on possibility space $(\Theta, \mathbb{P}(\Theta), Pos) \,\forall i \in \{1, 2, \ldots, n\}$. Then the primitive chance of rough fuzzy event characterized by $u_j(\xi) \leq 0, j = 1, 2, \ldots, m$, is a measurable map from $[0,1]$ to $[0,1]$ and is defined below in (17).
$$Ch\{u_j(\hat{\xi}) \leq 0, j = 1, 2, \ldots, m\}(\beta) = \sup\left\{\alpha \,\middle|\, Cr\left\{\theta \in \Theta \,\middle|\, Tr\left\{\begin{matrix}u_j(\hat{\xi}) \leq 0\\ j=1,2,\ldots,m\end{matrix}\right\} \geq \alpha\right\} \geq \beta\right\}, \quad (17)$$
where $\alpha$ and $\beta$ are the predetermined confidence limits.

## 3. Multi-objective optimization

A multi-objective optimization problem (MOP) optimizes a vector of objectives of conflicting nature.

MOP can be formulated as

$$\begin{aligned}
&\text{Minimize } Z(x) = \big(z_1(x), z_2(x), \ldots, z_m(x)\big) \in \Re^m \\
&\text{such that} \\
&g_i(x) \leq 0;\ i = 1, 2, \ldots, p \\
&h_j(x) = 0;\ j = p+1, p+2, \ldots, q \\
&x = (x_1, x_2, \ldots, x_n)^T \in \Re^n\ \forall x_l \in [x_l^L, x_l^U],
\end{aligned} \qquad (18)$$

where $x$ is a $n$ dimensional vector of decision variables, $f_k$ is the $k^{\text{th}}$ objective function, $g_i$ is the $i^{\text{th}}$ inequality constraint and $h_j$ is the $j^{\text{th}}$ equality constraint.

Due to the conflicting nature among the objectives in MOP there exist tradeoffs between them. Hence improving an objective diminishes the value of at least one of the remaining objectives. Therefore instead of a single optimal solution a MOP always generate an optimal set of solutions. These solutions cannot be compared with each other, i.e., one cannot recommend that a particular solution in the optimal set is better than the remaining solutions. The solutions of the optimal set are nondominated solutions in the sense that neither of the solutions can dominate each other.

Consider $r = (r_1, r_2, \ldots, r_n)^T$, $t = (t_1, t_2, \ldots, t_n)^T \in \Re^n$ are two solutions of an optimal set for the model (18), then $r$ is said to dominate $t$ (or $r \prec t$) if and only if,

    i. For any objective $z_k$, $z_k(r) \leq z_k(t), \forall\, k \in \{1, 2, \ldots, m\}$
    ii. There exists at least one objective $f_l$ such that $z_l(r) \leq z_l(t)$, where $l \in \{1, 2, \ldots, m\}$.

$r$ and $t$ are said to be nondominated to each other if neither $r \prec t$ nor $t \prec r$.

A nondominated solution $x^*$ is said to be Pareto optimal solution if there does not exists any other possible solution $x$ in the decision space such that $Z(x) \prec Z(x^*)$. The set of all Pareto optimal solutions generate the Pareto set (PS) and the set of all objective vectors in the objective space, corresponding to PS is known as Pareto front (PF).

## 4. Multi-objective model with rough fuzzy coefficients

Rough fuzzy multi-objective problem (RFMOP) deals with optimizing a vector of conflicting objectives subject to a set of constraints such that the coefficients are those associated with the decision variables of the problem are rough fuzzy in nature. A RFMOP model is depicted below:

$$\begin{cases} \text{Minimize } \big(f_1(\hat{\bar{\zeta}}_1^T, \mathcal{X}), f_2(\hat{\bar{\zeta}}_2^T, \mathcal{X}), \ldots, f_m(\hat{\bar{\zeta}}_m^T, \mathcal{X})\big) \\ \text{such that} \\ g_i(\hat{\bar{\zeta}}_i^T, \mathcal{X}) \leq 0;\ i = 1, 2, \ldots, p, \end{cases} \qquad (19)$$

where, $m$ is the number of objective functions, $\mathcal{X} = (x_1, x_2, \ldots, x_n)^T$ is the decision vector, $\hat{\bar{\zeta}}_l^T = (\hat{\bar{\zeta}}_{l1}, \hat{\bar{\zeta}}_{l2}, \ldots, \hat{\bar{\zeta}}_{ln})^T\ l = 1, 2, \ldots, m$; $\hat{\bar{\zeta}}_i^T = (\hat{\bar{\zeta}}_{i1}, \hat{\bar{\zeta}}_{i2}, \ldots, \hat{\bar{\zeta}}_{in})^T; i = 1, 2, \ldots, p$ are the rough fuzzy vectors.

From (17), the primitive chance of $g_i(\hat{\bar{\zeta}}_i^T, \mathcal{X}) \leq 0$ can be stated as:
$Ch\{g_i(\hat{\bar{\zeta}}_i^T, \mathcal{X}) \leq 0\}(\beta_i) \geq \alpha_i = Cr\{\theta \in \Theta | Tr\{g_i(\hat{\bar{\zeta}}_i^T, \mathcal{X}) \leq 0\} \geq \alpha_i\} \geq \beta_i,\ i = 1, 2, \ldots, p$
and
$Ch\{f_l(\hat{\bar{\zeta}}_l^T, \mathcal{X}) \leq \bar{\bar{f}}_l\}(\tau_l) \geq \rho_l = Cr\{\theta \in \Theta | Tr\{f_l(\hat{\bar{\zeta}}_l^T, \mathcal{X}) \leq \bar{\bar{f}}_l\} \geq \rho_l\} \geq \tau_l,\ l = 1, 2, \ldots, m.$

Furthermore, a decision vector $\mathcal{X}$ is said to be feasible to (19) if and only if the credibility measures of the fuzzy events, $\{\theta \in \Theta | Tr\{f_l(\hat{\bar{\zeta}}_l^T, \mathcal{X}) \leq \bar{\bar{f}}_l\} \geq \rho_l\}$ at confidence levels of at least $\tau_l$, are satisfied by a set of fuzzy event constraints, $\{\theta \in \Theta | Tr\{g_i(\hat{\bar{\zeta}}_i^T, \mathcal{X}) \leq 0\} \geq \alpha_i\}$ having at least credibility measures as $\beta_i$.

Hence, the chance constrained multi-objective programming model becomes:

$$\begin{cases} Minimize \ (\bar{\bar{f}}_1, \bar{\bar{f}}_2, \dots, \bar{\bar{f}}_m) \\ such\ that \\ Cr\{\theta \in \Theta | Tr\{f_l(\hat{\zeta}_l^T, \mathcal{X}) \leq \bar{\bar{f}}_l\} \geq \rho_l\} \geq \tau_l;\ l = 1,2, \dots, m \\ Cr\{\theta \in \Theta | Tr\{g_i(\hat{\zeta}_i^T, \mathcal{X}) \leq 0\} \geq \alpha_i\} \geq \beta_i;\ i = 1,2, \dots, p, \end{cases} \quad (20)$$

where $Cr\{\cdot\}$ and $Tr\{\cdot\}$ denote the credibility and trust measure of an event in $\{\cdot\}$ respectively. $\rho_l, \tau_l, \alpha_i$ and $\beta_i$ are the predetermined confidence levels with, $l = 1,2, \dots, m$ and $i = 1,2, \dots, p$.

Moreover, if $\mathcal{X}_1$ and $\mathcal{X}_2$ are two solutions of (20) at the predetermined confidence levels $\rho_l$ and $\tau_l$ then we can state that $\mathcal{X}_1$ dominates $\mathcal{X}_2$, i.e., $\mathcal{X}_1 \prec \mathcal{X}_2$ at a trust level of $\rho_l$ and credibility level of $\tau_l$ if and only if following two conditions are true.

i. $\{Cr\{\theta \in \Theta | Tr\{f_l(\hat{\zeta}_l^T, \mathcal{X}_1) \leq F_l(\mathcal{X}_1)\} \geq \rho_l\} \geq \tau_l\} \leq \{Cr\{\theta \in \Theta | Tr\{f_l(\hat{\zeta}_l^T, \mathcal{X}_2) \leq F_l(\mathcal{X}_1)\} \geq \rho_l\} \geq \tau_l\}, \forall l = 1,2, \dots, m$
ii. There exists at least one $j$ such that $F_j(\mathcal{X}_1) < F_j(\mathcal{X}_2)$, where $j \in \{1,2, \dots, m\}$

Again, $\mathcal{X}_1$ and $\mathcal{X}_2$ are said to be nondominated solutions at $\rho_l$ and $\tau_l$ if neither $\mathcal{X}_1 \prec \mathcal{X}_2$ nor $\mathcal{X}_2 \prec \mathcal{X}_1$ at trust level of $\rho_l$ and credibility level of $\tau_l$. For simplicity, we have expressed have expressed $F_l(\mathcal{X})$ as $\bar{\bar{f}}_l, \forall l = 1,2, \dots, m$.

**Remark 1:** A rough fuzzy vector $\hat{\zeta}_l$ is degenerated to a fuzzy vector, $\tilde{\zeta}_l$ if for $\rho_l > 0, Cr\{\theta \in \Theta | Tr\{f_l(\hat{\zeta}_l^T, \mathcal{X}) \leq \bar{\bar{f}}_l\} \geq \rho_l\} \geq \tau_l$ become equivalent to a fuzzy event $Cr\{\theta \in \Theta | f_l(\tilde{\zeta}_l^T, \mathcal{X})\} \geq \tau_l; l = 1,2, \dots, m$ which is a standard fuzzy chance constraint. Similarly, $Cr\{\theta \in \Theta | Tr\{g_i(\hat{\zeta}_i^T, \mathcal{X}) \leq 0\} \geq \alpha_i\} \geq \beta_i$ also degenerates to $Cr\{\theta \in \Theta | g_i(\tilde{\zeta}_i^T, \mathcal{X})\} \geq \beta_i$ at $\alpha_i > 0; i = 1,2, \dots, p$. Hence, the model (20) is transformed to an equivalent multi-objective fuzzy chance constrained programming model [25] as,

$$\begin{cases} Minimize\ (\bar{\bar{f}}_1, \bar{\bar{f}}_2, \dots, \bar{\bar{f}}_m) \\ such\ that \\ Cr\{\theta \in \Theta | f_l(\tilde{\zeta}_l^T, \mathcal{X}) \leq \bar{\bar{f}}_l\} \geq \tau_l;\ l = 1,2, \dots, m \\ Cr\{\theta \in \Theta | g_i(\tilde{\zeta}_i^T, \mathcal{X}) \leq 0\} \geq \beta_i;\ i = 1,2, \dots, p \end{cases}$$

**Remark 2:** A rough fuzzy vector $\hat{\zeta}_l$ is degenerated to a rough vector $\zeta_l$ if for $\tau_l > 0, Cr\{\theta \in \Theta | Tr\{f_l(\hat{\zeta}_l^T, \mathcal{X}) \leq \bar{\bar{f}}_l\} \geq \rho_l\} \geq \tau_l$ become equivalent to a rough event $Tr\{f_l(\zeta_l^T, \mathcal{X}) \leq \bar{\bar{f}}_l\} \geq \rho_l; l = 1,2, \dots, m$ which is a standard rough chance constraint. Similarly, $Cr\{\theta \in \Theta | Tr\{g_i(\hat{\zeta}_i^T, \mathcal{X}) \leq 0\} \geq \alpha_i\} \geq \beta_i$ degenerates to $Tr\{g_i(\zeta_i^T, \mathcal{X}) \leq 0\} \geq \alpha_i; i = 1,2, \dots, p$ at $\beta_i > 0$. Hence, the model (20) is transformed to an equivalent multi-objective rough chance constrained programming model [25] as,

$$\begin{cases} Minimize\ (\bar{\bar{f}}_1, \bar{\bar{f}}_2, \dots, \bar{\bar{f}}_m) \\ such\ that \\ Tr\{f_l(\zeta_l^T, \mathcal{X}) \leq \bar{\bar{f}}_l\} \geq \rho_l;\ l = 1,2, \dots, m \\ Tr\{g_i(\zeta_i^T, \mathcal{X}) \leq 0\} \geq \alpha_i;\ i = 1,2, \dots, p \end{cases}$$

**Lemma 4.1:** For any $\theta \in \Theta$, a rough fuzzy vector $\hat{\zeta}_{ij}\chi$, degenerated to a fuzzy vector $\tilde{\zeta}_i(\theta)^T\chi$, $\tilde{\zeta}_i(\theta) = [\tilde{\zeta}_{i1}(\theta), \tilde{\zeta}_{i2}(\theta), \dots, \tilde{\zeta}_{in}(\theta)]^T$ and $\chi = [x_1, x_2, \dots, x_n]^T$ are respectively a fuzzy vector and a decision vector.

**Proof:** Without loss of genareality, we have assumed that the fuzzy variable degenerated from a rough fuzzy variable for any $\theta \in \Theta$ becomes a triangular fuzzy variable. Let $\hat{\zeta}_{ij}$ be a rough fuzzy variable, then for any $\theta \in \Theta$ there exists a fuzzy variable $\tilde{\zeta}_{ij}(\theta) = (l_{ij}, r_{ij}, s_{ij})$ whose credibility measure is defined as,

$$Cr\big(\tilde{\zeta}_{ij}(\theta) \leq y\big) = \begin{cases} 1 & ; s_{ij} \leq y \\ \frac{y + s_{ij} - 2r_{ij}}{2(s_{ij} - r_{ij})} & ; r_{ij} \leq y \leq s_{ij} \\ \frac{y - l_{ij}}{2(r_{ij} - l_{ij})} & ; l_{ij} \leq y \leq r_{ij} \\ 0 & ; otherwise \end{cases}$$

and $Cr\big(\tilde{\zeta}_{ij}(\theta) \geq y\big) = \begin{cases} 1 & ; y \leq l_{ij} \\ \frac{2r_{ij} - l_{ij} - y}{2(r_{ij} - l_{ij})} & ; l_{ij} \leq y \leq r_{ij} \\ \frac{s_{ij} - y}{2(s_{ij} - r_{ij})} & ; r_{ij} \leq y \leq s_{ij} \\ 0 & ; otherwise \end{cases}$

From the extension principle it follows that

$$\tilde{\zeta}_i(\theta)^T \mathcal{X} = \sum_{j=1}^n \tilde{\zeta}_{ij}(\theta) \mathcal{X} = \left(\sum_{j=1}^n (l_{ij}\mathcal{X}_j, r_{ij}\mathcal{X}_j, s_{ij}\mathcal{X}_j)\right)$$

$$= \left(\sum_{j=1}^n l_{ij}\mathcal{X}_j, \sum_{j=1}^n r_{ij}\mathcal{X}_j, \sum_{j=1}^n s_{ij}\mathcal{X}_j\right) = (u, v, w),$$

where $u = \sum_{j=1}^n l_{ij}\mathcal{X}_j$, $v = \sum_{j=1}^n r_{ij}\mathcal{X}_j$ and $w = \sum_{j=1}^n s_{ij}\mathcal{X}_j$.

**Lemma 4.2:** For any $\theta \in \Theta$ a rough fuzzy vector $\hat{\zeta}_{ij}\mathcal{X}$, degenerated to a rough vector $\zeta_i(\theta)^T \mathcal{X}$, $\zeta_i(\theta)^T = [\zeta_{i1}(\theta), \zeta_{i2}(\theta), ..., \zeta_{in}(\theta)]^T$ and $\mathcal{X} = [x_1, x_2, ..., x_n]^T$ are respectively a rough vector and a decision vector.

**Proof:** For a given $\theta \in \Theta$, a rough fuzzy variable $\hat{\zeta}_{ij}$ when degenerated to a rough variable $\zeta_{ij}(\theta) = [u_{ij}, v_{ij}][s_{ij}, t_{ij}]$ with its trust distribution,

$$Tr\{\zeta_{ij}(\theta) \leq z\} = \begin{cases} 0 & ; z \leq s_{ij} \\ \frac{z - s_{ij}}{2(t_{ij} - s_{ij})} & ; s_{ij} \leq z \leq u_{ij} \\ \frac{1}{2}\left(\frac{z - u_{ij}}{v_{ij} - u_{ij}} + \frac{z - s_{ij}}{t_{ij} - s_{ij}}\right) & ; u_{ij} \leq z \leq v_{ij} \\ \frac{1}{2}\left(\frac{z - s_{ij}}{t_{ij} - s_{ij}} + 1\right) & ; v_{ij} \leq z \leq t_{ij} \\ 1 & ; z \geq t_{ij} \end{cases}$$

and $Tr\{\zeta_{ij}(\theta) \geq z\} = \begin{cases} 0 & ; z \geq q \\ \frac{t_{ij} - z}{2(t_{ij} - s_{ij})} & ; v_{ij} \leq z \leq t_{ij} \\ \frac{1}{2}\left(\frac{t_{ij} - z}{t_{ij} - s_{ij}} + \frac{v_{ij} - z}{v_{ij} - u_{ij}}\right) & ; u_{ij} \leq z \leq v_{ij} \\ \frac{1}{2}\left(\frac{t_{ij} - z}{t_{ij} - s_{ij}} + 1\right) & ; s_{ij} \leq z \leq u_{ij} \\ 1 & ; z \leq s_{ij} \end{cases}$

Then $\zeta_i(\theta)^T \mathcal{X} = \sum_{j=1}^n \zeta_{ij}(\theta) \mathcal{X} = \sum_{j=1}^n \big([u_{ij}\mathcal{X}_j, v_{ij}\mathcal{X}_j], [s_{ij}\mathcal{X}_j, t_{ij}\mathcal{X}_j]\big)$

$$= \left(\left[\sum_{j=1}^{n} u_{ij}\mathcal{X}_j, \sum_{j=1}^{n} v_{ij}\mathcal{X}_j\right], \left[\sum_{j=1}^{n} s_{ij}\mathcal{X}_j, \sum_{j=1}^{n} t_{ij}\mathcal{X}_j\right]\right) = [m,n][p,q], \text{ where } p \leq m \leq n \leq q \text{ such}$$
that $m = \sum_{j=1}^{n} u_{ij}\mathcal{X}_j, n = \sum_{j=1}^{n} v_{ij}\mathcal{X}_j, p = \sum_{j=1}^{n} s_{ij}\mathcal{X}_j$ and $q = \sum_{j=1}^{n} t_{ij}\mathcal{X}_j$.

**Theorem 4.1:** For given confidence levels $\alpha_i, \beta_i \in [0,1]$,

$Cr\{\theta \in \Theta | Tr\{\zeta_i(\theta)^T \mathcal{X} \leq \bar{\bar{f}}_i\} \geq \alpha_i\} \geq \beta_i$ implies $Cr\{\theta \in \Theta | p + 2\alpha_i Q_1 \leq \bar{\bar{f}}_i\} \geq \beta_i$, where $Q_1 = (q - p)$ and $p \leq \bar{\bar{f}}_i \leq m$.

**Proof:** From the result of Lemma 4.2, $\zeta_i(\theta)^T \mathcal{X}$ is known to be a rough variable.

Therefore the trust distribution function of $\zeta_i(\theta)^T \mathcal{X}$ is as given below when $p \leq \bar{\bar{f}}_i \leq m$.

$$Tr\{\zeta_i(\theta)^T \mathcal{X} \leq \bar{\bar{f}}_i\} = \begin{cases} 0 & ; \bar{\bar{f}}_i \leq p \\ 0.5 \frac{\bar{\bar{f}}_i - p}{(q-p)} & ; p \leq \bar{\bar{f}}_i \leq m \\ 0.5 \left(\frac{\bar{\bar{f}}_i - m}{n-m} + \frac{\bar{\bar{f}}_i - p}{q-p}\right) & ; m \leq \bar{\bar{f}}_i \leq n \\ 0.5 \left(\frac{\bar{\bar{f}}_i - p}{q-p} + 1\right) & ; n \leq \bar{\bar{f}}_i \leq q \\ 1 & ; \bar{\bar{f}}_i \geq q \end{cases}$$

$$\Rightarrow \alpha_i \leq \begin{cases} 0 & ; \bar{\bar{f}}_i \leq p \\ 0.5 \frac{\bar{\bar{f}}_i - p}{(q-p)} & ; p \leq \bar{\bar{f}}_i \leq m \\ 0.5 \left(\frac{\bar{\bar{f}}_i - m}{n-m} + \frac{\bar{\bar{f}}_i - p}{q-p}\right) & ; m \leq \bar{\bar{f}}_i \leq n \\ 0.5 \left(\frac{\bar{\bar{f}}_i - p}{q-p} + 1\right) & ; n \leq \bar{\bar{f}}_i \leq q \\ 1 & ; \bar{\bar{f}}_i \geq q \end{cases}$$

Here, $\alpha_i$ is the predetermined confidence level of the trust distribution for the rough variable $\zeta_i(\theta)^T \mathcal{X}$.

It immediately follows that

$\bar{\bar{f}}_i \geq p + 2\alpha_i(q - p)$ if $p \leq \bar{\bar{f}}_i \leq m$
$\bar{\bar{f}}_i \geq \frac{mq + pn + 2\alpha_i(n-m)(q-p) - 2pm}{n+q-(p+m)}$ if $m \leq \bar{\bar{f}}_i \leq n$
and
$\bar{\bar{f}}_i \geq 2p - q + 2\alpha_i(q - p)$, if $n \leq \bar{\bar{f}}_i \leq q$

Therefore considering all the cases together, we get
$Cr\{\theta \in \Theta | Tr\{\zeta_i(\theta)^T \mathcal{X} \leq \bar{\bar{f}}_i\} \geq \alpha_i\} \geq \beta_i$

$$\Rightarrow Cr\left\{\theta \in \Theta | \bar{\bar{f}}_i \geq \begin{cases} p + 2\alpha_i(q-p) & ; p \leq \bar{\bar{f}}_i \leq m \\ \frac{mq+pn+2\alpha_i(n-m)(q-p)-2pm}{n+q-(p+m)} & ; m \leq \bar{\bar{f}}_i \leq n \\ 2p - q + 2\alpha_i(q-p) & ; n \leq \bar{\bar{f}}_i \leq q \\ q & ; q \leq \bar{\bar{f}}_i \end{cases}\right\} \geq \beta_i$$

$$\Rightarrow Cr\left\{\theta \in \Theta | \bar{\bar{f}}_i \geq \begin{cases} p + 2\alpha_i Q_1 & ; p \leq \bar{\bar{f}}_i \leq m \\ \frac{mq+pn+2\alpha_i Q_2 - 2pm}{n+q-(p+m)} & ; m \leq \bar{\bar{f}}_i \leq n \\ 2p - q + 2\alpha_i Q_1 & ; n \leq \bar{\bar{f}}_i \leq q \\ q & ; q \leq \bar{\bar{f}}_i \end{cases}\right\} \geq \beta_i$$

where, $Q_1 = (q - p)$ and $Q_2 = (n - m)(q - p)$.

Hence, it follows from the above result that

$Cr\{\theta \in \Theta | Tr\{\zeta_i(\theta)^T \mathcal{X} \leq \bar{\bar{f}}_i\} \geq \alpha_i\} \geq \beta_i \Rightarrow Cr\{\theta \in \Theta | p + 2\alpha_i Q_1 \leq \bar{\bar{f}}_i\} \geq \beta_i$ where, $p \leq \bar{\bar{f}}_i \leq m$ and $Q_1 = (q - p)$.

**Theorem 4.2:** For given confidence levels $\alpha_i, \beta_i \in [0,1]$, when $p \leq \bar{\bar{f}}_i \leq m$

$$Cr\{\theta \in \Theta | Tr\{\zeta_i(\theta)^T \mathcal{X} \leq \bar{\bar{f}}_i\} \geq \alpha_i\} \geq \beta_i$$

$$\Rightarrow \bar{\bar{f}}_i \geq \begin{cases} w + 2\alpha_i w_{Q_1} & ; w \leq \bar{\bar{f}}_i - 2\alpha_i Q_1 \\ 2(v + \alpha_i v_{Q_1}) - (w + 2\alpha_i w_{Q_1}) + 2\beta_i\left((w + 2\alpha_i w_{Q_1}) - (v + 2\alpha_i v_{Q_1})\right) & ; v \leq \bar{\bar{f}}_i - 2\alpha_i Q_1 \leq w \\ u + 2\alpha_i u_{Q_1} + 2\beta_i\left((v + 2\alpha_i v_{Q_1}) - (u + 2\alpha_i u_{Q_1})\right) & ; u \leq \bar{\bar{f}}_i - 2\alpha_i Q_1 \leq v \end{cases}$$

Here, $Q_1 = (q - p)$ is a TFV of the form, $(u_{Q_1}, v_{Q_1}, w_{Q_1})$.

**Proof:** From the result of Theorem 4.1, $Cr\{\theta \in \Theta | Tr\{\zeta_i(\theta)^T \mathcal{X} \leq \bar{\bar{f}}_i\} \geq \alpha_i\} \geq \beta_i \Rightarrow Cr\{\theta \in \Theta | p + 2\alpha_i Q_1 \leq \bar{\bar{f}}_i\} \geq \beta_i$ when $p \leq \bar{\bar{f}}_i \leq m$ and $Q_1 = (q - p)$. Now $p$ and $Q_1$ become fuzzy variables of the form $\tilde{\zeta}_i(\theta)^T \mathcal{X}$. The fuzzy variable in the rest of the paper is assumed to be the triangular fuzzy variable of the form, $(u, v, w)$, cf. Lemma 4.1 and thus,

$Cr\{\theta \in \Theta | p + 2\alpha_i Q_1 \leq \bar{\bar{f}}_i\} \geq \beta_i \Rightarrow Cr\{\theta \in \Theta | p \leq \bar{\bar{f}}_i - 2\alpha_i Q_1\} \geq \beta_i$ when $p \leq \bar{\bar{f}}_i - 2\alpha_i Q_1 \leq m$

$$\Rightarrow \beta_i \leq \begin{cases} 1 & ; w \leq \bar{\bar{f}}_i - 2\alpha_i Q_1 \\ \frac{\bar{\bar{f}}_i - 2\alpha_i Q_1 + w - 2v}{2(w-v)} & ; v \leq \bar{\bar{f}}_i - 2\alpha_i Q_1 \leq w \\ \frac{\bar{\bar{f}}_i - 2\alpha_i Q_1 - u}{2(v-u)} & ; u \leq \bar{\bar{f}}_i - 2\alpha_i Q_1 \leq v \\ 0 & ; otherwise \end{cases}$$

$$\Rightarrow \bar{\bar{f}}_i \geq \begin{cases} w + 2\alpha_i Q_1 & ; w \leq \bar{\bar{f}}_i - 2\alpha_i Q_1 \\ 2v - w + 2\alpha_i Q_1 + 2\beta_i(w - v) & ; v \leq \bar{\bar{f}}_i - 2\alpha_i Q_1 \leq w \\ u + 2\alpha_i Q_1 + 2\beta_i(v - u) & ; u \leq \bar{\bar{f}}_i - 2\alpha_i Q_1 \leq v \end{cases}$$

$$\Rightarrow \bar{\bar{f}}_i \geq \begin{cases} w + 2\alpha_i(q - p) & ; w \leq \bar{\bar{f}}_i - 2\alpha_i Q_1 \\ 2v - w + 2\alpha_i(q - p) + 2\beta_i(w - v) & ; v \leq \bar{\bar{f}}_i - 2\alpha_i Q_1 \leq w \\ u + 2\alpha_i(q - p) + 2\beta_i(v - u) & ; u \leq \bar{\bar{f}}_i - 2\alpha_i Q_1 \leq v \end{cases} \Rightarrow \bar{\bar{f}}_i \geq$$

$$\begin{cases} w + 2\alpha_i w_{Q_1} & ; w \leq \bar{\bar{f}}_i - 2\alpha_i Q_1 \\ 2(v + \alpha_i v_{Q_1}) - (w + 2\alpha_i w_{Q_1}) + 2\beta_i\left((w + 2\alpha_i w_{Q_1}) - (v + 2\alpha_i v_{Q_1})\right) & ; v \leq \bar{\bar{f}}_i - 2\alpha_i Q_1 \leq w \\ u + 2\alpha_i u_{Q_1} + 2\beta_i\left((v + 2\alpha_i v_{Q_1}) - (u + 2\alpha_i u_{Q_1})\right) & ; u \leq \bar{\bar{f}}_i - 2\alpha_i Q_1 \leq v \end{cases}$$

Here $Q_1$ is a TFV expressed as $(u_{Q_1}, v_{Q_1}, w_{Q_1})$

Similarly, we can also have,

$Cr\{\theta \in \Theta | Tr\{\zeta_i(\theta)^T \mathcal{X} \leq \bar{\bar{f}}_i\} \geq \alpha_i\} \geq \beta_i \Rightarrow Cr\left\{\theta \in \Theta | \frac{mq+pn+2\alpha_i Q_2 - 2pm}{n+q-(p+m)} \leq \bar{\bar{f}}_i\right\} \geq \beta_i$ when $m \leq \bar{\bar{f}}_i \leq n$ and $Q_2 = (n-m)(q-p)$ and $Cr\{\theta \in \Theta | Tr\{\zeta_i(\theta)^T \mathcal{X} \leq \bar{\bar{f}}_i\} \geq \alpha_i\} \geq \beta_i \Rightarrow Cr\{\theta \in \Theta | 2p - q + 2\alpha_i Q_1 \leq \bar{\bar{f}}_i\} \geq \beta_i$ when $n \leq \bar{\bar{f}}_i \leq q$ and $Q_1 = (q-p)$

Thenceforth,

$Cr\left\{\theta \in \Theta | \frac{mq+pn+2\alpha_i Q_2 - 2pm}{n+q-(p+m)} \leq \bar{\bar{f}}_i\right\} \geq \beta_i$

$\Rightarrow Cr\left\{\theta \in \Theta | \frac{mq+pn-2pm}{n+q-(p+m)} \leq \bar{\bar{f}}_i - 2\alpha_i U\right\} \geq \beta_i$

$\Rightarrow \bar{\bar{f}}_i \geq \begin{cases} w + 2\alpha_i U & ; w \leq \bar{\bar{f}}_i - 2\alpha_i U \\ 2v - w + 2\alpha_i U + 2\beta_i(w-v) & ; v \leq \bar{\bar{f}}_i - 2\alpha_i U \leq w \\ u + 2\alpha_i U + 2\beta_i(v-u) & ; u \leq \bar{\bar{f}}_i - 2\alpha_i U \leq v \end{cases}$

$\Rightarrow \bar{\bar{f}}_i \geq \begin{cases} w + \frac{2\alpha_i(n-m)(q-p)}{n+q-(p+m)} & ; w \leq \bar{\bar{f}}_i - 2\alpha_i U \\ 2v - w + \frac{2\alpha_i(n-m)(q-p)}{n+q-(p+m)} + 2\beta_i(w-v) & ; v \leq \bar{\bar{f}}_i - 2\alpha_i U \leq w \\ u + \frac{2\alpha_i(n-m)(q-p)}{n+q-(p+m)} + 2\beta_i(v-u) & ; u \leq \bar{\bar{f}}_i - 2\alpha_i U \leq v, \end{cases}$

$\Rightarrow \bar{\bar{f}}_i \geq \begin{cases} w + 2\alpha_i w_U & ; w \leq \bar{\bar{f}}_i - 2\alpha_i U \\ 2(v + \alpha_i v_U) - (w + 2\alpha_i w_U) + 2\beta_i\big((w + 2\alpha_i w_U) - (v + 2\alpha_i v_U)\big) & ; v \leq \bar{\bar{f}}_i - 2\alpha_i U \leq w \\ u + 2\alpha_i u_U + 2\beta_i\big((v + 2\alpha_i v_U) - (u + 2\alpha_i u_U)\big) & ; u \leq \bar{\bar{f}}_i - 2\alpha_i U \leq v \end{cases}$

where $U = \frac{Q_2}{n+q-(p+m)}$ is a triangular fuzzy variable of the form, $(u_U, v_U, w_U)$

Also,

$Cr\{\theta \in \Theta | 2p - q + 2\alpha_i Q_1 \leq \bar{\bar{f}}_i\} \geq \beta_i$

$\Rightarrow Cr\{\theta \in \Theta | 2p - q \leq \bar{\bar{f}}_i - 2\alpha_i Q_1\} \geq \beta_i$

$\Rightarrow \bar{\bar{f}}_i \geq \begin{cases} w + 2\alpha_i Q_1 & ; w \leq \bar{\bar{f}}_i - 2\alpha_i Q_1 \\ 2v - w + 2\alpha_i Q_1 + 2\beta_i(w-v) & ; v \leq \bar{\bar{f}}_i - 2\alpha_i Q_1 \leq w \\ u + 2\alpha_i Q_1 + 2\beta_i(v-u) & ; u \leq \bar{\bar{f}}_i - 2\alpha_i Q_1 \leq v \end{cases}$

Since, $p, q$ and $Q_1$ are TFVs then, $2p - q + 2\alpha_i Q_1$ will also become a TFV, which immediately follows,

$\Rightarrow \bar{\bar{f}}_i \geq \begin{cases} w + 2\alpha_i(q-p) & ; w \leq \bar{\bar{f}}_i - 2\alpha_i Q_1 \\ 2v - w + 2\alpha_i(q-p) + 2\beta_i(w-v) & ; v \leq \bar{\bar{f}}_i - 2\alpha_i Q_1 \leq w \\ u + 2\alpha_i(q-p) + 2\beta_i(v-u) & ; u \leq \bar{\bar{f}}_i - 2\alpha_i Q_1 \leq v \end{cases}$

$$\Rightarrow \bar{\bar{f}}_i \geq \begin{cases} w + 2\alpha_i w_{Q_1} & ; w \leq \bar{\bar{f}}_i - 2\alpha_i Q_1 \\ 2(v + \alpha_i v_{Q_1}) - (w + 2\alpha_i w_{Q_1}) + 2\beta_i \left( (w + 2\alpha_i w_{Q_1}) - (v + 2\alpha_i v_{Q_1}) \right) & ; v \leq \bar{\bar{f}}_i - 2\alpha_i Q_1 \leq w \\ u + 2\alpha_i u_{Q_1} + 2\beta_i \left( (v + 2\alpha_i v_{Q_1}) - (u + 2\alpha_i u_{Q_1}) \right) & ; u \leq \bar{\bar{f}}_i - 2\alpha_i Q_1 \leq v \end{cases}$$

**Theorem 4.3:** For any given confidence level $\alpha_i \in [0,1]$

(i) If $0 \leq \beta_i \leq 0.5$, then $Cr\{\theta \in \Theta | p \leq \bar{\bar{f}}_i - 2\alpha_i Q_1\} \geq \beta_i$ if and only if,
$u + 2\alpha_i u_{Q_1} + 2\beta_i \{(v + 2\alpha_i v_{Q_1}) - (u + 2\alpha_i u_{Q_1})\} \leq \bar{\bar{f}}_i$ for $u \leq \bar{\bar{f}}_i - 2\alpha_i Q_1 \leq v$

(ii) If $0.5 \leq \beta_k \leq 1$, then $Cr\{\theta \in \Theta | p \leq \bar{\bar{f}}_i - 2\alpha_i Q_1\} \geq \beta_i$ if and only if,
$2(v + \alpha_i v_{Q_1}) - (w + 2\alpha_i w_{Q_1}) + 2\beta_i \left( (w + 2\alpha_i w_{Q_1}) - (v + 2\alpha_i v_{Q_1}) \right) \leq \bar{\bar{f}}_i$
for $v \leq \bar{\bar{f}}_i - 2\alpha_i Q_1 \leq w$, where $Q_1 = (q - p)$ is of the form $(u_{Q_1}, v_{Q_1}, w_{Q_1})$.

**Proof:**

i. By Lemma 2.1 and by Theorem 4.2, if $0 \leq \beta_i \leq 0.5$ and $Cr\{\theta \in \Theta | p + Q_1 \leq \bar{\bar{f}}_i\} \geq \beta_i$ then
$Cr\{\theta \in \Theta | p \leq \bar{\bar{f}}_i - Q_1\} \geq \beta_i \Rightarrow u \leq \bar{\bar{f}}_i - Q_1 \leq v$

Henceforth by Theorem 4.2,
$Cr\{\theta \in \Theta | p \leq \bar{\bar{f}}_i - 2\alpha_i Q_1\} \geq \beta_i \Rightarrow (u + 2\alpha_i Q_1 + 2\beta_i(v - u)) \leq \bar{\bar{f}}_i$
$\Rightarrow \bar{\bar{f}}_i \geq (u(1 - 2\beta_i) + 2\alpha_i(q - p) + 2\beta_i v)$
$\Rightarrow \bar{\bar{f}}_i \geq \{u + 2\alpha_i u_{Q_1} + 2\beta_i \{(v + 2\alpha_i v_{Q_1}) - (u + 2\alpha_i u_{Q_1})\}\}$

Conversly, if $u \leq \bar{\bar{f}}_i - 2\alpha_i Q_1 \leq v$ then $0 \leq Cr\{\theta \in \Theta | p \leq \bar{\bar{f}}_i - 2\alpha_i Q_1\} \geq \beta_i \geq \frac{1}{2}$

$\Rightarrow \frac{1}{2} \leq \beta_i \leq \left( \frac{\bar{\bar{f}}_i - 2\alpha_i Q_1 - u}{2(v-u)} \right) \Rightarrow Cr\{\theta \in \Theta | p \leq \bar{\bar{f}}_i - 2\alpha_i Q_1\} \geq \beta_i$

ii. If $0.5 \leq \beta_i \leq 1$ and $Cr\{\theta \in \Theta | p \leq \bar{\bar{f}}_i - 2\alpha_i Q_1\} \geq \beta_i$ then by Lemma 2.1 and Theorem 4.2,
$Cr\{\theta \in \Theta | p \leq \bar{\bar{f}}_i - 2\alpha_i Q_1\} \geq \beta_i \Rightarrow v \leq (\bar{\bar{f}}_i - 2\alpha_i Q_1) \leq w$
Then, by Theorem 4.2,
$Cr\{\theta \in \Theta | p \leq \bar{\bar{f}}_i - 2\alpha_i Q_1\} \geq \beta_i \Rightarrow v \leq (\bar{\bar{f}}_i - 2\alpha_i Q_1) \leq w$

Again by Theorem 4.2,
$Cr\{\theta \in \Theta | p \leq \bar{\bar{f}}_i - Q_1\} \geq \beta_i \Rightarrow \bar{\bar{f}}_i \geq (2v - w + 2\alpha_i Q_1 + 2\beta_i(w - v))$
$\Rightarrow (2(1 - \beta_i)v + (2\beta_i - 1)w + 2\alpha_i(q - p)) \leq \bar{\bar{f}}_i$
$\Rightarrow \{2(v + \alpha_i v_{Q_1}) - (w + 2\alpha_i w_{Q_1}) + 2\beta_i \left( (w + 2\alpha_i w_{Q_1}) - (v + 2\alpha_i v_{Q_1}) \right)\} \leq \bar{\bar{f}}_i$

Conversly, if $v \leq \bar{\bar{f}}_i - 2\alpha_i Q_1 \leq w$, then $\frac{1}{2} \leq Cr\{\theta \in \Theta | p \leq \bar{\bar{f}}_i - 2\alpha_i Q_1\} = 1 \geq \beta_i$

$\Rightarrow \frac{1}{2} \leq \beta_i \leq \left( \frac{\bar{\bar{f}}_i - 2\alpha_i Q_1 + w - 2v}{2(w-v)} \right) = 1$
$\Rightarrow Cr\{\theta \in \Theta | p - 2\alpha_i Q_1 \leq \bar{\bar{f}}_i\} \geq \beta_i$

Similarly when $0 \leq \beta_i \leq 0.5$ then for $u \leq (\bar{\bar{f}}_i - 2\alpha_i U) \leq v$,
$$Cr\left\{\theta \in \Theta | \frac{mq+pn-2pm}{n+q-(p+m)} \leq \bar{\bar{f}}_i - 2\alpha_i U\right\} \geq \beta_i$$
$$\Rightarrow (2(1-\beta_i)v + (2\beta_i - 1)w + 2\alpha_i U) \leq \bar{\bar{f}}_i$$
$$\Rightarrow u + 2\alpha_i u_U + 2\beta_i\{(v + 2\alpha_i v_U) - (u + 2\alpha_i u_U)\}$$

and for $u \leq \bar{\bar{f}}_i - Q_1 \leq v$,
$$Cr\{\theta \in \Theta | 2p - q \leq \bar{\bar{f}}_i - 2\alpha_i Q_1\} \geq \beta_i$$
$$\Rightarrow (u(1 - 2\beta_i) + 2\beta_i v + 2\alpha_i Q_1) \leq \bar{\bar{f}}_i$$

Furthermore when $0.5 \leq \beta_i \leq 1$ then for $u \leq (\bar{\bar{f}}_i - 2\alpha_i U) \leq v$,
$$Cr\left\{\theta \in \Theta | \frac{mq+pn-2pm}{n+q-(p+m)} \leq \bar{\bar{f}}_i - 2\alpha_i U\right\} \geq \beta_i$$
$$\Rightarrow \bar{\bar{f}}_i \geq ((1-2\beta_i)u + 2\alpha_i U + 2\beta_i v)$$
$$\Rightarrow 2(v + \alpha_i v_U) - (w + 2\alpha_i w_U) + 2\beta_i((w + 2\alpha_i w_U) - (v + 2\alpha_i v_U))$$

and for $v \leq (\bar{\bar{f}}_i - 2\alpha_i Q_1) \leq w$,
$$Cr\{\theta \in \Theta | 2p - q \leq \bar{\bar{f}}_i - 2\alpha_i Q_1\} \geq \beta_i$$
$$\Rightarrow (2(1-\beta_i)v + (2\beta_i - 1)w + 2\alpha_i Q_1) \leq \bar{\bar{f}}_i,$$

where $Q_1 = (q - p)$ and $U = \frac{(n-m)(q-p)}{n+q-(p+m)}$.

## 5. Quadratic Minimum Spanning Tree Problem with Rough Fuzzy coefficients

In this section, a bi-objective $(\alpha, \beta)$ rough fuzzy quadratic minimum spanning tree problem (b-$(\alpha, \beta)$RFQMSTP) is modeled for a WCN $G = (V_G, E_G)$, where $V_G$ is the vertex set and $E_G$ is the edge set of $G$. Two types of weights are associated with $G$: linear and quadratic. The linear weights are associated with the edges, while the quadratic weights are allied with the interaction or interfacing of a pair of edges in $G$. In this paper, we have considered the weight as cost. The linear cost $\xi_e$ of an edge $e(\in E_G)$ is the cost incurred while traversing $e$ and the quadratic cost $\zeta_{ee'}$ implies the interaction/interfacing cost between any pair of edges $e$ and $e'$.

For example, in a utility or transportation network, the linear weight resembles the expenditure related to manufacturing and installation of gas, water pipes, cost incurred for buying construction materials, labour charges, etc. Whereas, the quadratic weights represent the interaction/interfacing costs for installing valves or bending losses in T-joints in a utility transmission network, or turn delays in transportation network [34]. The interfacing costs mostly occur for adjacent edges in these cases. The interfacing cost may also associate with any pair of edges in case of fibre-optic networks [34].

In such real networks, there always exists uncertainty in the values of the cost parameters. Different experts may give different opinions about the cost and the opinion of the experts may be represented in range or intervals. In such cases the uncertainty associated with the information may be imprecise and as well as vague. Such uncertainty can be expressed by a rough fuzzy variable [5, 6].

In this article, we have considered the linear costs and the quadratic costs as rough fuzzy variables [5, 6] and used the chance constrained programming to model b-$(\alpha, \beta)$ RFQMSTP and transformed the model to its crisp equivalent. Henceforth, we first rewrite the chance constrained multi-objective programming model, represented earlier in (20), as follows in (21) for bi-objective case.

$$\begin{cases} \text{Minimize } (\bar{\bar{f}}_1, \bar{\bar{f}}_2) \\ \text{such that} \\ Cr\{\theta \in \Theta | Tr\{f_1(\xi_s(\theta), x_s) \leq \bar{\bar{f}}_1\} \geq \alpha_1\} \geq \beta_1 \\ Cr\{\theta \in \Theta | Tr\{f_2(\zeta_{st}(\theta), x_s, x_t) \leq \bar{\bar{f}}_2\} \geq \alpha_2\} \geq \beta_2 \\ \sum_{s=1}^{|E_G|} x_s = |V_G| - 1 \\ \sum_{s \in E_\kappa} x_s \leq |\kappa| - 1, \ \kappa \subset V_G, |\kappa| \geq 3 \\ \forall \ x_s, x_t \in \{0,1\}; \ x_s \neq x_t \\ \alpha_1, \alpha_2, \beta_1, \beta_2 \in [0,1] \end{cases} \quad (21)$$

Here, $\alpha_1, \alpha_2, \beta_1$ and $\beta_2$ are predetermined confidence levels and $|V_G|$ and $|E_G|$ are respectively the order and size of $G$. Problem (21) can again be expressed as

$$\begin{cases} \text{Minimize } (\bar{\bar{f}}_1, \bar{\bar{f}}_2) \\ \text{such that} \\ Cr\left\{\theta \in \Theta | Tr\left\{\sum_{s=1}^{|E_G|} \xi_s(\theta) x_s \leq \bar{\bar{f}}_1\right\} \geq \alpha_1\right\} \geq \beta_1 \\ Cr\left\{\theta \in \Theta | Tr\left\{\sum_{s=1}^{|E_G|} \sum_{t=1}^{|E_G|} \zeta_{st}(\theta) x_s x_t \leq \bar{\bar{f}}_2\right\} \geq \alpha_2\right\} \geq \beta_2 \\ \sum_{s=1}^{|E_G|} x_s = |V_G| - 1 \\ \sum_{s \in E_\kappa} x_s \leq |\kappa| - 1, \ \kappa \subset V_G, |\kappa| \geq 3 \\ \forall \ x_s, x_t \in \{0,1\}; \ x_s \neq x_t \\ \alpha_1, \alpha_2, \beta_1, \beta_2 \in [0,1] \end{cases} \quad (22)$$

Let $L(\theta) = \sum_{s=1}^{|E_G|} \xi_s(\theta) x_s$ and $Q(\theta) = \sum_{s=1}^{|E_G|} \sum_{t=1}^{|E_G|} \zeta_{st}(\theta) x_s x_t$ then $L(\theta)$ and $Q(\theta)$ become rough variables similar to $[m,n][p,q]$ by Lemma 4.2. Hence $L(\theta)$ and $Q(\theta)$ becomes $[m_1, n_1][p_1, q_1]$ and $[m_2, n_2][p_2, q_2]$ respectively. Here according to Theorem 4.1, we have modeled problem (23) by considering the cases, $p_1 \leq L(\theta) \leq m_1$ and $p_2 \leq Q(\theta) \leq m_2$.

Then problem defined in (22) becomes equivalent to

$$\begin{cases} \text{Minimize } (\bar{\bar{f}}_1, \bar{\bar{f}}_2) \\ \text{such that} \\ Cr\{\theta \in \Theta | p_1 + 2\alpha_1(q_1 - p_1) \leq \bar{\bar{f}}_1\} \geq \beta_1 \\ Cr\{\theta \in \Theta | p_2 + 2\alpha_2(q_2 - p_2) \leq \bar{\bar{f}}_2\} \geq \beta_2 \\ \sum_{s=1}^{|E_G|} x_s = |V_G| - 1 \\ \sum_{s \in E_\kappa} x_s \leq |\kappa| - 1, \ \kappa \subset V_G, |\kappa| \geq 3 \\ \forall \ x_s, x_t \in \{0,1\}; \ x_s \neq x_t \\ \alpha_1, \alpha_2, \beta_1, \beta_2 \in [0,1] \end{cases} \quad (23)$$

Again, $\tilde{L}(\theta) = p_1 + 2\alpha_1(q_1 - p_1)$ and $\tilde{Q}(\theta) = p_2 + 2\alpha_2(q_2 - p_2)$ become the TFVs such that $p_1 = (u_1, v_1, w_1)$, $(q_1 - p_1) = (u_L, v_L, w_L)$ and $p_2 = (u_2, v_2, w_2)$ and $(q_2 - p_2) = (u_Q, v_Q, w_Q)$. Then by Theorem 4.2 and Theorem 4.3, model (23) is now equivalent to (24) and (25) as follows.

$$\begin{cases} \bar{\bar{f}}_1 = \text{Minimize } u_1 + 2\alpha_1 u_L + 2\beta_1\left((v_1 + 2\alpha_1 v_L) - (u_1 + 2\alpha_1 u_L)\right) \\ \bar{\bar{f}}_2 = \text{Minimize } u_2 + 2\alpha_2 u_Q + 2\beta_1\left((v_2 + 2\alpha_2 v_Q) - (u_2 + 2\alpha_2 u_Q)\right) \\ such\ that \\ \sum_{s=1}^{|E_G|} x_s = |V_G| - 1 \\ \sum_{s \in E_\kappa} x_s \leq |\kappa| - 1,\ \kappa \subset V_G, |\kappa| \geq 3 \\ \forall\ x_s, x_t \in \{0,1\};\ x_s \neq x_t \\ \alpha_1, \alpha_2 \in [0,1];\ \beta_1, \beta_2 \in [0,0.5] \end{cases} \quad (24)$$

and

$$\begin{cases} \bar{\bar{f}}_1 = \text{Minimize } 2(v_1 + \alpha_1 v_L) - (w_1 + 2\alpha_2 w_L) + 2\beta_1\left((w_1 + 2\alpha_1 w_L) - (v_1 + 2\alpha_1 v_L)\right) \\ \bar{\bar{f}}_2 = \text{Minimize } 2(v_2 + \alpha_2 v_Q) - (w_2 + 2\alpha_2 w_Q) + 2\beta_2\left((w_2 + 2\alpha_2 w_Q) - (v_2 + 2\alpha_2 v_Q)\right) \\ such\ that \\ \sum_{s=1}^{|E_G|} x_s = |V_G| - 1 \\ \sum_{s \in E_\kappa} x_s \leq |\kappa| - 1,\ \kappa \subset V_G, |\kappa| \geq 3 \\ \forall\ x_s, x_t \in \{0,1\};\ x_s \neq x_t \\ \alpha_1, \alpha_2 \in [0,1];\ \beta_1, \beta_2 \in [0.5,1] \end{cases} \quad (25)$$

## 6. Methodologies for solving bi-objective $(\alpha, \beta)$ rough fuzzy quadratic minimum spanning tree problem

A bi-objective quadratic minimum spanning tree problem has been considered with rough fuzzy coefficients. The model of the problem is formulated using chance-constrained programming technique (cf. problem (21)) which is transformed to the crisp equivalent form (cf. problem (25)). This crisp equivalent model of problem (21) is then solved with epsilon-constraint approach, which is a classical approach to solve a MOP. Furthermore, the same model is also solved with two MOEAs, namely, NSGA-II [27] and MOCHC [23].

### 6.1. Epsilon-constraint

Haimes et al. [39] proposed epsilon-constraint technique to solve the multi-objective optimization problems. In this approach a MOP is re-modelled to optimize a single objective function and remaining objectives are restricted to a user defined values. This method is also proved to be efficient enough to find solutions of the problems having non convex objective space. The corresponding re-formulated model of problem (18) that can be solved by epsilon-constraint method is as follows.

$$\begin{aligned} & \text{Minimize } F(x) = f_r(x); \\ & such\ that \\ & \quad f_s(x) \leq \varepsilon_s;\ s = 1, 2, \dots, M\ \text{and}\ s \neq r \\ & \quad g_i(x) \leq 0;\ i = 1, 2, \dots, p \\ & \quad h_j(x) = 0;\ j = p+1, p+2, \dots, q \\ & \quad x = (x_1, x_2, \dots, x_n)^T \in \Re^n\ \forall x_l \in [x_l^L, x_l^U], \end{aligned} \quad (26)$$

where x is a vector of decision variables, $f_k$ is the $k^{th}$ objective function, $g_i$ is $i^{th}$ inequality constraint and $h_j$ is $j^{th}$ equality constraint. $\varepsilon_s$ is the user defined value and set to upper bound of $f_s(x)$ whose proximity is not necessarily zero.

### 6.2. Multi-objective Genetic Algorithm (MOGA)

MOGA, proposed by Fonseca and Flaming [10], has drawn colossal attention in the recent past for solving complex multi-objective optimization problems. MOGAs assert nondominated solutions in a generational population and at the same time promoting diversity among the solutions. Various MOGAs exist in the literature such as Horn et al.[21], Srinivas and Deb [32], Zitzler and Thiele [13], Deb et al. [23], Zang and Li [40], Coello et al. [8, 9], Konak et al. [3], Nag et al. [28], etc. A MOGA generates a set of solutions in each iteration, which create a nondominated front in the objective space. The Pareto front (PF) is the ideal one corresponding to first nondominated front of the solutions. In this paper, we consider two MOGAs: NSGA-II and MOCHC to solve bi-objective $(\alpha, \beta)$ rough fuzzy quadratic minimum spanning tree problem.

**6.2.1. Nondominated Sorting Genetic Algorithm II (NSGA-II)**

The nondominated sorting genetic algorithm II (NSGA-II) [27] is a stochastic population based elitist multi-objective genetic algorithm. Along with elite preservation strategy, NSGA-II also emphasizes on diversity preservation. NSGA-II assures that every generational population retains certain proportion of the elite chromosomes from the population of previous generation. In a generation, the parent of size $N$ produces same number of offspring using selection, crossover and mutation operators. The offspring are then intermingled with the parent to form a total of $2N$ solutions. Among these, $N$ number of best solutions are finally considered for next generation using crowded-comparison operator [32]. The crowded-comparison operator uses two metrics: nondomination rank ($i_{rank}$) and crowding distance ($i_{distance}$). Solutions having lower $i_{rank}$ are preferred over higher $i_{rank}$. If $i_{rank}$ of the solutions are identical then the solutions with greater $i_{distance}$, *i.e.*, solutions from less congested regions are selected. Thus NSGA-II promotes both the elitism and diversity in every generation. This mechanism is applied in every generation until the termination criteria is achieved. The termination criteria can be a preset threshold of either number of fitness evaluations or computer clock time or diversity in the generational population. The computational complexity of NSGA-II is $O(MN^2)$, where $M$ and $N$ are the number of objectives and the population size respectively. The working principle of NSGA-II when used to solve problems (24) and (25) is shown in Fig 1.

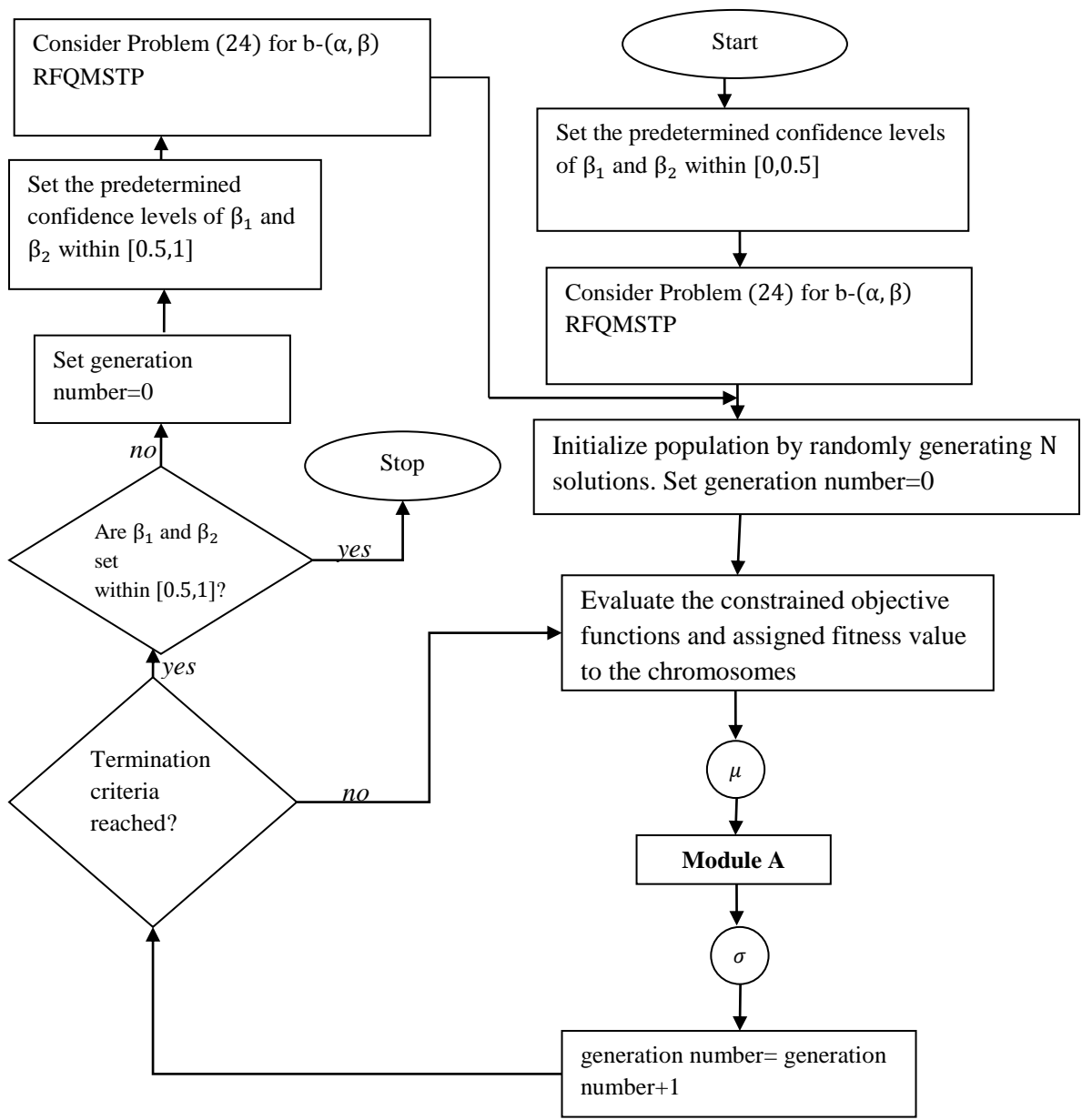

Fig. 1: Flowchart of NSGA-II for solving b-$(\alpha, \beta)$ RFQMSTP

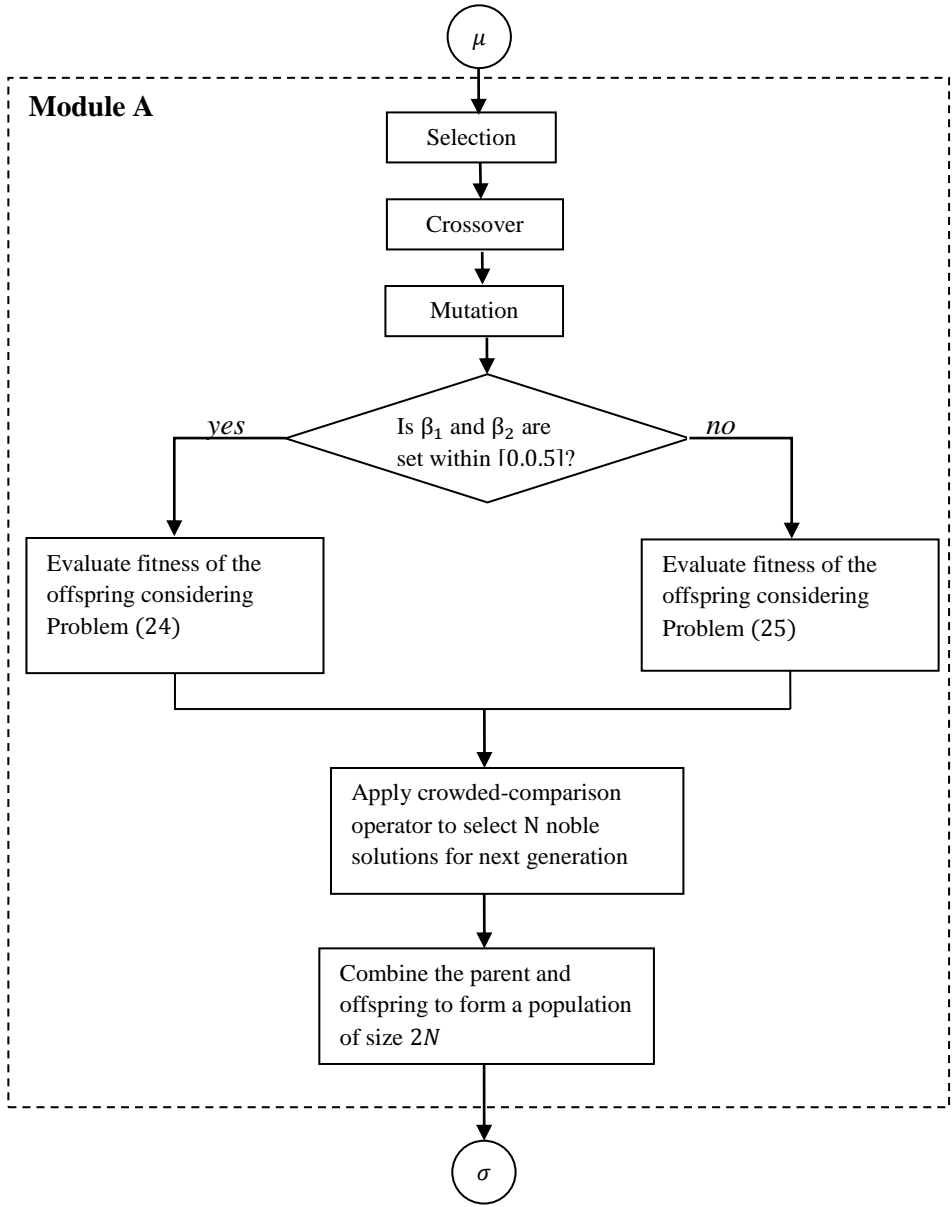

Fig. 2: A modular component of flowchart in Fig. 1

### 6.2.2. Multi-objective Cross Generational Elitist Selection, Heterogeneous Recombination and Cataclysmic Mutation (MOCHC)

The multi-objective CHC algorithm proposed by Nebro et al. [23] is the multi-objective version of Eshelman's CHC algorithm [29]. MOCHC favours elitism highly, where the offspring have to compete with the parents for survival. Along with the cross-generational elitist strategy, MOCHC also accentuates on recombination operator, i.e., crossover for diversity preservation while mutation is only used in the *restart* phase [23] of the algorithm in order to avoid the pre-mature convergence of the algorithm. MOCHC can be interpreted as a multi-objective genetic algorithm where offsprings are produced by a recombination operator known as half uniform crossover (HUX) [29] in subsequent generations without any mutation operator until the *restart* [23] phase of the algorithm is reached.

The MOCHC starts by randomly generating a population of size $N$. In a generation, the selections of parents for recombination mechanism are done randomly in a way to implement incest prevention, i.e., prohibition of mating of similar solutions in order to maintain diversity in the population. Only those parents which are not similar can participate in a recombination mechanism, i.e., if the Hamming distance between the parents are greater than a threshold limit [23] only then the parents are considered for recombination. In MOCHC the recombination mechanism which is applied on the selected parents is commonly known as HUX [29]. The HUX [29] mechanism replicates the common information of the parents to their respective offspring and then exchanges half of the diverging information from each parent to the produced offspring so that the Hamming distance between the parents and the offspring is maximum. Here, the objective of this mechanism is to introduce diversity in the population since there is no other way to add diversity in the population of a generation due to the absence of mutation operator. Once the offspring are generated, they are then combined with their parent population among which the best (elite) $N$ nondominated solutions are selected by crowded-comparison operator [32]. If $N$ nondominated solutions for the two consecutive generations remain unaltered then the threshold value also decreases by one. In this way, as the generation increases the threshold limit also progressively decreases since the population become more homogeneous (i.e., with the increase in generation the chances of generating new offspring will also reduce reasonably) with generations. When the threshold value [23] becomes 0 it is assumed that no new diversity is introduced in the population and the population become complete homogeneous. The reason for such complete homogeneity is the premature convergence of the population due to elitism and absence of mutation. At this stage, in order to add new diversity in the population the *restart* phase of MOCHC begins. In the *restart* phase, best 5% non dominated solutions of the present population are preserved which are selected by crowded-comparison [32] operator from the current population and the rest of the solutions are cataclysmically mutated using bit-flip mutation with a very high mutation probability (0.35) as suggested by Eshelman [29] and Nebro et al. [23]. High mutation probability is used to introduce new diversity in the population. Once the *restart* phase is over, then same recombination and elitism mechanism are applied on the population and this entire process continues until the termination criteria is satisfied. Fig. 3 diagrammatically explains the working principle of MOCHC when used to solve problems: (24) and (25).

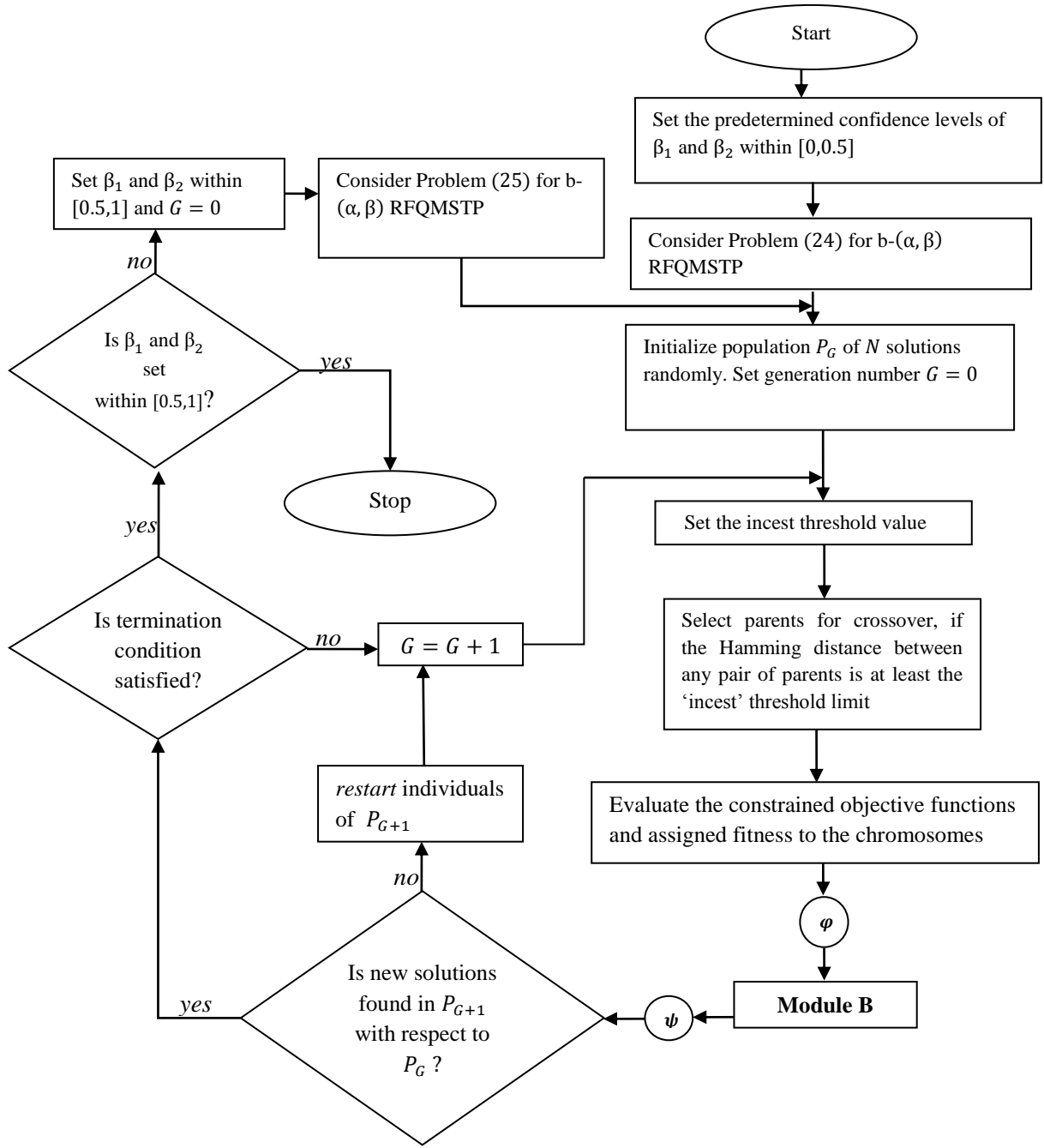

Fig. 3: Flowchart of MOCHC for solving b-$(\alpha, \beta)$ RFQMSTP

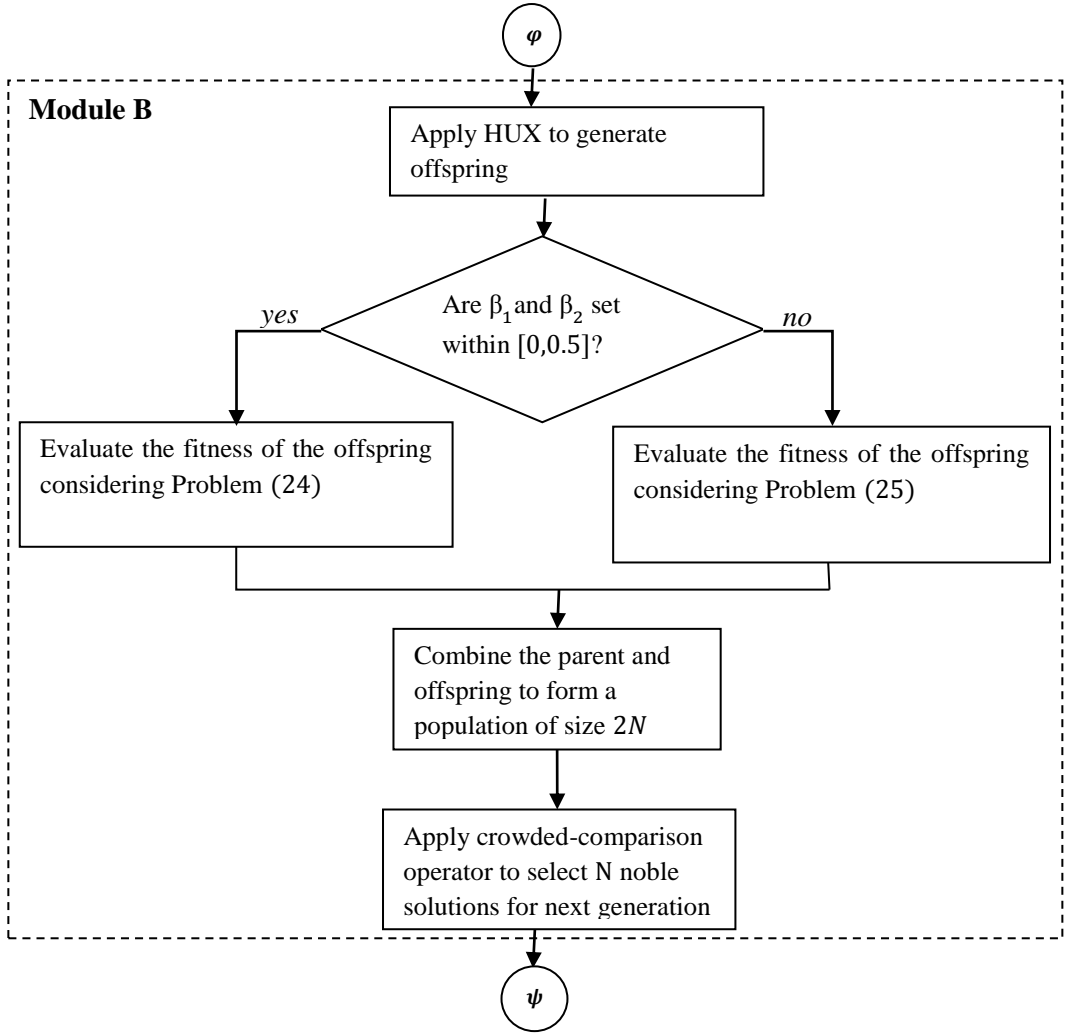

Fig. 4: A modular component of flowchart in Fig. 3

## 7. Results and Discussion:

We have consider a weighted connected network (WCN) as shown in Fig. 5, $G = (V_G, E_G)$ where $V_G$ is the vertex set and $E_G$ is the edge set of $G$. The corresponding linear and quadratic weights of $G$ are expressed as rough fuzzy variables.

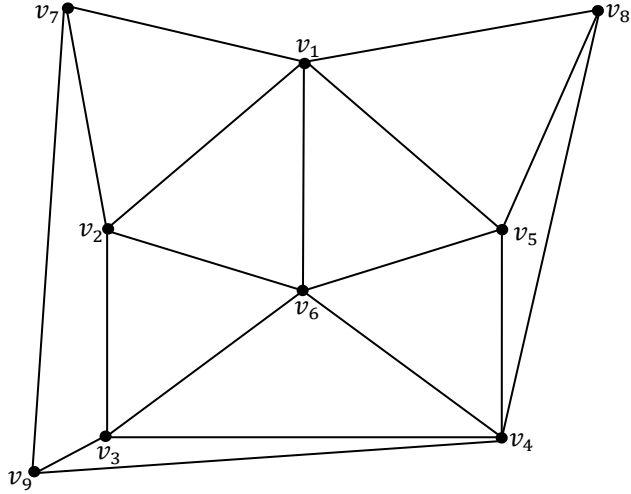

Fig. 5: A weighted connected graph $G$

Table 1: Rough fuzzy linear weights for $G$

| Edges | Rough fuzzy linear weights | Edges | Rough fuzzy linear weights |
|---|---|---|---|
| $e_{12}$ | $[\xi, \xi+2][\xi-1, \xi+3]$ $\xi = (9, 11.5, 12.7)$ | $e_{56}$ | $[\xi, \xi+0.5][\xi-0.7, \xi+0.9]$ $\xi = (9.4, 11.6, 12.8)$ |
| $e_{15}$ | $[\xi, \xi+1][\xi-1, \xi+2]$ $\xi = (11, 13, 15)$ | $e_{17}$ | $[\xi, \xi+0.9][\xi-1.2, \xi+1.8]$ $\xi = (11.5, 12.9, 13.8)$ |
| $e_{16}$ | $[\xi, \xi+1.5][\xi-0.7, \xi+1.9]$ $\xi = (9.2, 10.5, 12.4)$ | $e_{18}$ | $[\xi, \xi+0.5][\xi-0.6, \xi+0.95]$ $\xi = (8.9, 10.2, 12)$ |
| $e_{23}$ | $[\xi, \xi+2][\xi-2, \xi+3]$ $\xi = (10.6, 14.1, 17.2)$ | $e_{27}$ | $[\xi, \xi+1][\xi-1.1, \xi+1.5]$ $\xi = (10, 12, 13.5)$ |
| $e_{26}$ | $[\xi, \xi+1][\xi-2, \xi+2]$ $\xi = (8.3, 10, 12.4)$ | $e_{39}$ | $[\xi, \xi+0.4][\xi-0.5, \xi+0.8]$ $\xi = (9.8, 11.2, 12.4)$ |
| $e_{34}$ | $[\xi, \xi+2][\xi-2, \xi+2.5]$ $\xi = (9.8, 11.7, 14)$ | $e_{58}$ | $[\xi, \xi+1.2][\xi-1.2, \xi+1.4]$ $\xi = (11.9, 12.8, 14.5)$ |
| $e_{36}$ | $[\xi, \xi+1.5][\xi-0.5, \xi+2]$ $\xi = (11, 14, 15)$ | $e_{79}$ | $[\xi, \xi+1.1][\xi-0.9, \xi+1.3]$ $\xi = (11.5, 12, 13.5)$ |
| $e_{45}$ | $[\xi, \xi+1][\xi-0.5, \xi+1.8]$ $\xi = (12, 14.7, 17.1)$ | $e_{48}$ | $[\xi, \xi+0.6][\xi-0.4, \xi+0.67]$ $\xi = (9, 10.2, 12.9)$ |
| $e_{46}$ | $[\xi, \xi+2.5][\xi-1.5, \xi+2.7]$ $\xi = (10, 12, 14)$ | $e_{49}$ | $[\xi, \xi+0.78][\xi-0.8, \xi+0.8]$ $\xi = (10.2, 11.6, 12.4)$ |

The linear and quadratic weights of $G$ are presented in Table 1 and Table 2 respectively. Table 2 displays quadratic weights of $G$. The remaining quadratic weights of $G$ which are not present in Table 2 are assumed to be zero for simplicity. We consider the linear and quadratic weights of $G$ as two different objectives to be minimized. This bi-objective model formulation (b-$(\alpha, \beta)$ RFQMSTP) is done by chance-constrained programming technique, presented earlier in problem (21), given the rough fuzzy weights and pre determined confidence levels $\alpha$ and $\beta$.

Table 2: Rough fuzzy quadratic weights for $G$

| Edges | Rough fuzzy quadratic weights | Edges | Rough fuzzy quadratic weights |
|---|---|---|---|
| $e_{12}e_{15}$ | $[\zeta, \zeta + 0.7][\xi - 0.2, \zeta + 1]$ $\zeta = (9, 11, 13)$ | $e_{56}e_{36}$ | $[\zeta, \zeta + 0.9][\zeta - 0.7, \zeta + 1.2]$ $\zeta = (10, 12.5, 14)$ |
| $e_{12}e_{16}$ | $[\zeta, \zeta + 1.2][\zeta - 0.5, \zeta + 1.4]$ $\zeta = (8.5, 10.2, 11.5)$ | $e_{12}e_{56}$ | $[\zeta, \zeta + 0.7][\zeta - 0.5, \zeta + 1.0]$ $\zeta = (9, 11, 13)$ |
| $e_{15}e_{16}$ | $[\zeta, \zeta + 0.4][\zeta - 0.2, \zeta + 0.8]$ $\zeta = (9.5, 10.7, 11.6)$ | $e_{15}e_{58}$ | $[\zeta, \zeta + 0.5][\zeta - 0.5, \zeta + 0.9]$ $\zeta = (10.9, 11.5, 12.1)$ |
| $e_{15}e_{56}$ | $[\zeta, \zeta + 1][\zeta - 0.8, \zeta + 1.2]$ $\zeta = (9.8, 10.8, 11.5)$ | $e_{46}e_{48}$ | $[\zeta, \zeta + 0.9][\zeta - 0.4, \zeta + 1.2]$ $\zeta = (10.8, 12.2, 13.5)$ |
| $e_{26}e_{36}$ | $[\zeta, \zeta + 1.2][\zeta - 0.3, \zeta + 1.5]$ $\xi = (10.2, 12.2, 13.4)$ | $e_{15}e_{49}$ | $[\zeta, \zeta + 0.7][\zeta - 0.3, \zeta + 0.8]$ $\zeta = (11.2, 12, 13.9)$ |
| $e_{23}e_{26}$ | $[\zeta, \zeta + 0.9][\zeta - 0.5, \zeta + 1.6]$ $\zeta = (10.4, 11.7, 12.9)$ | $e_{17}e_{48}$ | $[\zeta, \zeta + 1.2][\zeta - 1.1, \zeta + 1.1]$ $\zeta = (10, 11.5, 12.5)$ |
| $e_{12}e_{26}$ | $[\zeta, \zeta + 1.2][\zeta - 0.2, \zeta + 1.5]$ $\zeta = (8.6, 9.2, 9.8)$ | $e_{17}e_{18}$ | $[\zeta, \zeta + 0.7][\zeta - 0.6, \zeta + 0.8]$ $\zeta = (8, 10, 12)$ |
| $e_{26}e_{46}$ | $[\zeta, \zeta + 0.9][\zeta - 0.3, \zeta + 1.2]$ $\zeta = (8.1, 9.8, 10.9)$ | $e_{18}e_{79}$ | $[\zeta, \zeta + 1.5][\zeta - 1.3, \zeta + 1.7]$ $\zeta = (8.9, 10, 11.2)$ |
| $e_{23}e_{34}$ | $[\zeta, \zeta + 2][\zeta - 0.9, \zeta + 2.6]$ $\zeta = (10.8, 12.9, 14.2)$ | $e_{27}e_{39}$ | $[\zeta, \zeta + 1.2][\zeta - 1.2, \zeta + 1.4]$ $\zeta = (11.2, 12.3, 13.4)$ |
| $e_{23}e_{46}$ | $[\zeta, \zeta + 2.1][\zeta - 0.2, \zeta + 2.4]$ $\zeta = (11, 13, 14)$ | $e_{27}e_{58}$ | $[\zeta, \zeta + 1.6][\zeta - 1.5, \zeta + 1.7]$ $\zeta = (10, 12.2, 14)$ |
| $e_{23}e_{56}$ | $[\zeta, \zeta + 1][\zeta - 1.2, \zeta + 1.5]$ $\zeta = (10.7, 11.7, 12.6)$ | $e_{39}e_{58}$ | $[\zeta, \zeta + 1][\zeta - 1, \zeta + 1.2]$ $\zeta = (9, 11, 13)$ |
| $e_{34}e_{56}$ | $[\zeta, \zeta + 1][\zeta - 1, \zeta + 2]$ $\zeta = (8, 10, 12)$ | $e_{27}e_{56}$ | $[\zeta, \zeta + 0.2][\zeta - 0.8, \zeta + 1]$ $\zeta = (12.1, 12.6, 12.9)$ |
| $e_{16}e_{34}$ | $[\zeta, \zeta + 0.9][\zeta - 1, \zeta + 1.9]$ $\zeta = (7.8, 11, 12)$ | $e_{27}e_{46}$ | $[\zeta, \zeta + 2][\zeta - 1, \zeta + 3]$ $\zeta = (10, 11, 12)$ |
| $e_{34}e_{46}$ | $[\zeta, \zeta + 1][\zeta - 1, \zeta + 2]$ $\zeta = (9, 11, 13)$ | $e_{39}e_{58}$ | $[\zeta, \zeta + 1][\zeta - 1, \zeta + 2]$ $\zeta = (10.5, 11.5, 12.9)$ |
| $e_{48}e_{79}$ | $[\zeta, \zeta + 1.5][\zeta - 1.2, \zeta + 1.9]$ $\zeta = (9.8, 11.2, 12.4)$ | $e_{26}e_{39}$ | $[\zeta, \zeta + 1.2][\zeta - 2, \zeta + 2.4]$ $\zeta = (11, 12.6, 13.5)$ |
| $e_{45}e_{36}$ | $[\zeta, \zeta + 2][\zeta - 1, \zeta + 3]$ $\zeta = (12, 14, 15)$ | $e_{46}e_{58}$ | $[\zeta, \zeta + 1.8][\zeta - 1.3, \zeta + 2.9]$ $\zeta = (12, 13.6, 14.5)$ |
| $e_{45}e_{46}$ | $[\zeta, \zeta + 0.7][\zeta - 0.5, \zeta + 0.9]$ $\zeta = (11.6, 12.7, 13.8)$ | $e_{49}e_{56}$ | $[\zeta, \zeta + 1][\zeta - 1.4, \zeta + 1.9]$ $\zeta = (8.9, 10.9, 12.5)$ |
| $e_{16}e_{56}$ | $[\zeta, \zeta + 1.4][\zeta - 0.1, \zeta + 1.5]$ $\zeta = (9.5, 10.5, 11.8)$ | $e_{45}e_{56}$ | $[\zeta, \zeta + 0.5][\zeta - 0.4, \zeta + 0.6]$ $\zeta = (10.8, 11.8, 12.9)$ |

We apply the methodologies, discussed in Section 5, to determine the QMST and its corresponding nondominated solutions corresponding to the graph $G$. We have solved problem (24) and problem (25) which are essentially the crisp equivalents of problem (21) for $G$ by epsilon-constraint method using a standard optimization solver LINGO. Table 3 shows the solutions in term of compromised decision vector $x$ and its corresponding objective values when (24) and (25) are solved by epsilon-constraint method. The confidence levels $\alpha$ and $\beta$ for both the objectives of problem (24) and problem (25) are set as $[\alpha = 0.9, \beta = 0.4]$ and $[\alpha = 0.9, \beta = 0.8]$ respectively. Fig. 6 and Fig. 7 portrays the QMSTs of $G$ corresponding to the optimal decision vectors of Table 3. The edges of these QMSTs are represented as dotted lines.

Table 3: Solutions of problem (24) and problem (25) with respect to $G$

| Confidence Levels $(\alpha, \beta)$ | Objective values, $\langle \bar{\bar{f}}_1, \bar{\bar{f}}_2 \rangle$ | Optimal decision vector, $x$ |
|---|---|---|
| $\alpha = 0.9, \beta = 0.4$ | $\langle 128.5600, 11.9400 \rangle$ | $(1\ 0\ 0\ 0\ 1\ 0\ 0\ 0\ 1\ 1\ 1\ 0\ 1\ 1\ 0\ 1\ 0\ 0)^T$ |
| $\alpha = 0.9, \beta = 0.8$ | $\langle 129.5760, 16.4400 \rangle$ | $(0\ 0\ 1\ 0\ 0\ 0\ 1\ 0\ 0\ 1\ 0\ 1\ 1\ 0\ 1\ 0\ 1\ 1)^T$ |

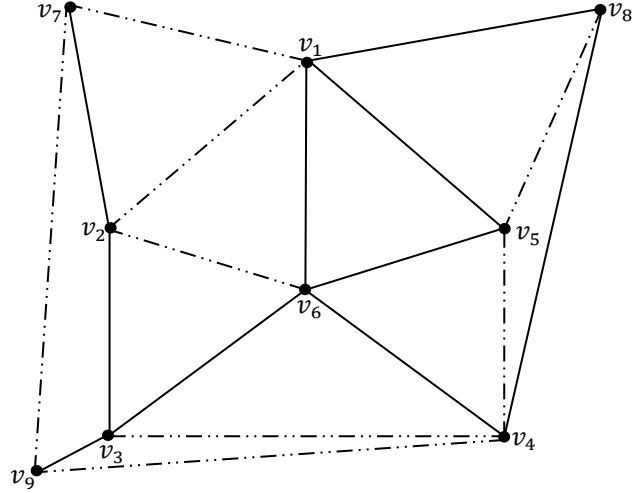

Fig. 6: A rough fuzzy $(\mathbf{0.9, 0.4})$ quadratic minimum spanning tree of $\boldsymbol{G}$

For both the problems, defined in (24) and (25), two different combination of values for $(\alpha, \beta)$ pair generate two distinct compromised binary decision vectors. A decision variable in a decision vector corresponds to an edge in $G$. If a decision variable is set to 1 the corresponding edge is considered in QMST of $G$, otherwise the edge is not included in the QMST. Fig. 6 shows the QMST corresponding to the compromised decision vector in the first row in Table 3. The QMST for $G$ displayed in Fig. 6 contains only the edges $e_{12}, e_{17}, e_{26}, e_{34}, e_{45}, e_{49}, e_{58}$ and $e_{79}$ for which the value of the corresponding decision variables are set to 1. Similarly, Fig. 7 displays the QMST corresponding to the compromised decision vector in the second row of Table 3. It shows that $\alpha$ and $\beta$ plays a significant role to generate the results for the crisp equivalent models of b-$(\alpha, \beta)$RFQMSTP.

The nondominated solutions obtained by solving problem (24) with confidence levels, $\alpha = 0.9$, $\beta = 0.4$ after 25,000 function evaluations of NSGA-II [27] and MOCHC [23] are displayed in Fig. 8 (a) and Fig 8(b) respectively. Again the corresponding figures of Fig. 9 (a) and Fig. 9 (b) depict the nondominated solutions after 25,000 function evaluations of NSGA-II [27] and MOCHC [23] for problem (25) with confidence levels, $\alpha = 0.9$, $\beta = 0.8$. Considering figures, 8(a)-(b) and 9 (a)-(b) we observe that the solution obtained using epsilon-constraint method is also nondominated compare to the solutions of NSGA-II [27] and MOCHC [23].

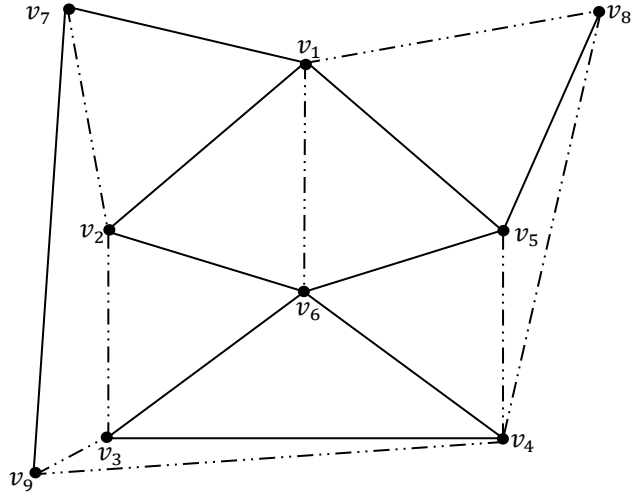

Fig. 7: A rough fuzzy ($0.9, 0.8$) quadratic minimum spanning tree of $G$

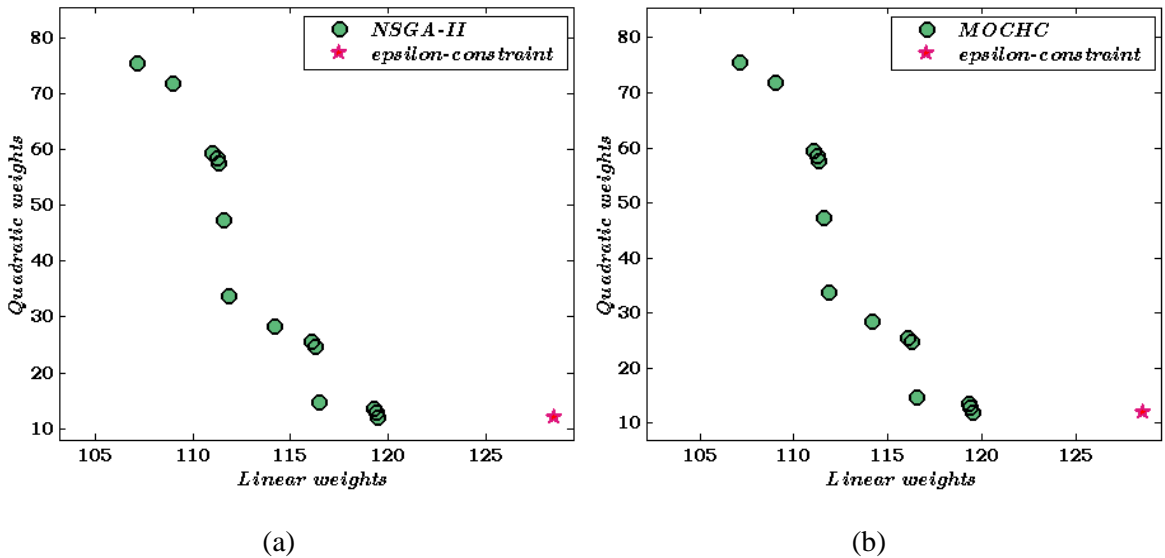

(a)                  (b)

Fig. 8: (a) Nondominated solutions of problem (24) for predetermined confidence levels of $\alpha = 0.9$ and $\beta = 0.4$ when solved by NSGA-II and epsilon-constraint method. (b) Nondominated solutions of problem (24) for predetermined confidence levels of $\alpha = 0.9$ and $\beta = 0.4$ when solved by MOCHC and epsilon-constraint method.

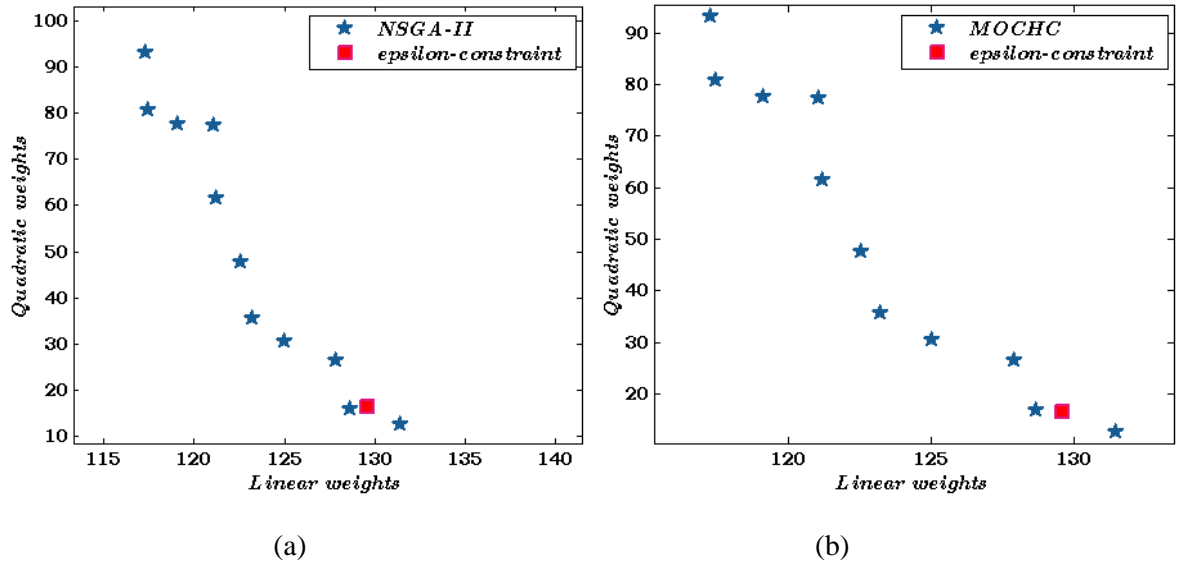

(a)                          (b)

Fig. 9: (a) Nondominated solutions of problem (25) for predetermined confidence levels of $\alpha = 0.9$ and $\beta = 0.8$ when solved by NSGA-II and epsilon-constraint method. (b) Nondominated solutions of problem (25) for predetermined confidence levels of $\alpha = 0.9$ and $\beta = 0.8$ when solved by MOCHC and epsilon-constraint method.

Table 4: Solutions of problem (24) and problem (25) for different confidence levels of $\alpha$ and $\beta$ for $G$

| $Cr = \beta$ | $Tr = \alpha$ | $\langle \bar{\bar{f}}_1, \bar{\bar{f}}_2 \rangle$ | Optimal decision vector, $x$ |
|---|---|---|---|
| 0.1 | 0.2 | ⟨87.4, 9.06⟩ | $(1\ 0\ 1\ 0\ 0\ 0\ 1\ 0\ 0\ 1\ 1\ 0\ 1\ 1\ 0\ 1\ 0\ 0)^T$ |
| | 0.4 | ⟨91.28, 11.36⟩ | $(0\ 0\ 1\ 0\ 1\ 0\ 1\ 0\ 0\ 1\ 1\ 0\ 1\ 1\ 0\ 1\ 0\ 0)^T$ |
| | 0.6 | ⟨100.00, 14.00⟩ | $(1\ 0\ 0\ 0\ 0\ 1\ 0\ 0\ 1\ 1\ 1\ 0\ 1\ 1\ 0\ 1\ 0\ 0)^T$ |
| | 0.8 | ⟨104.56, 10.86⟩ | $(1\ 0\ 1\ 0\ 0\ 0\ 1\ 0\ 0\ 1\ 1\ 0\ 1\ 1\ 0\ 1\ 0\ 0)^T$ |
| 0.3 | 0.2 | ⟨94.32, 9.66⟩ | $(1\ 0\ 0\ 0\ 1\ 0\ 0\ 0\ 1\ 1\ 1\ 0\ 1\ 1\ 0\ 1\ 0\ 0)^T$ |
| | 0.4 | ⟨98.39, 20.12⟩ | $(0\ 1\ 0\ 0\ 1\ 0\ 0\ 0\ 1\ 1\ 0\ 0\ 0\ 1\ 1\ 1\ 1\ 0)^T$ |
| | 0.6 | ⟨104.76, 10.14⟩ | $(1\ 0\ 0\ 0\ 1\ 0\ 0\ 0\ 1\ 1\ 1\ 0\ 1\ 1\ 0\ 1\ 0\ 0)^T$ |
| | 0.8 | ⟨109.97, 22.66⟩ | $(0\ 1\ 0\ 0\ 1\ 0\ 0\ 0\ 1\ 1\ 0\ 0\ 0\ 1\ 1\ 1\ 1\ 0)^T$ |
| 0.5 | 0.2 | ⟨104.00, 10.34⟩ | $(1\ 0\ 0\ 0\ 1\ 0\ 0\ 0\ 1\ 1\ 1\ 0\ 1\ 1\ 0\ 1\ 0\ 0)^T$ |
| | 0.4 | ⟨109.2, 10.58⟩ | $(0\ 0\ 0\ 0\ 1\ 0\ 1\ 0\ 1\ 1\ 1\ 0\ 1\ 1\ 0\ 1\ 0\ 0)^T$ |
| | 0.6 | ⟨108.14, 10.82⟩ | $(1\ 0\ 0\ 0\ 1\ 0\ 0\ 0\ 1\ 1\ 0\ 1\ 1\ 0\ 1\ 0\ 1\ 0)^T$ |
| | 0.8 | ⟨119.23, 24.42⟩ | $(0\ 1\ 0\ 0\ 1\ 0\ 0\ 0\ 1\ 1\ 0\ 0\ 0\ 1\ 1\ 1\ 1\ 0)^T$ |
| 0.7 | 0.2 | ⟨104.2480, 14.08⟩ | $(0\ 0\ 1\ 0\ 0\ 0\ 1\ 0\ 0\ 1\ 0\ 1\ 1\ 0\ 1\ 0\ 1\ 1)^T$ |
| | 0.4 | ⟨116.1180, 27.24⟩ | $(0\ 0\ 1\ 0\ 0\ 0\ 1\ 0\ 0\ 1\ 0\ 1\ 1\ 1\ 0\ 1\ 1\ 0)^T$ |
| | 0.6 | ⟨114.2680, 14.26⟩ | $(0\ 0\ 1\ 0\ 1\ 0\ 0\ 0\ 1\ 1\ 0\ 1\ 1\ 0\ 1\ 0\ 1\ 0)^T$ |
| | 0.8 | ⟨122.7880, 16.26⟩ | $(1\ 0\ 0\ 0\ 0\ 1\ 0\ 0\ 1\ 0\ 1\ 1\ 1\ 0\ 1\ 0\ 1\ 0)^T$ |
| 0.9 | 0.2 | ⟨111.0960, 13.30⟩ | $(0\ 0\ 1\ 0\ 1\ 0\ 0\ 0\ 1\ 1\ 0\ 1\ 1\ 0\ 1\ 0\ 1\ 0)^T$ |
| | 0.4 | ⟨128.56, 11.94⟩ | $(1\ 0\ 0\ 0\ 1\ 0\ 0\ 0\ 1\ 1\ 1\ 0\ 1\ 1\ 0\ 1\ 0\ 0)^T$ |
| | 0.6 | ⟨127.91, 26.66⟩ | $(0\ 0\ 1\ 0\ 0\ 0\ 1\ 0\ 0\ 1\ 0\ 0\ 1\ 1\ 1\ 1\ 0\ 1)^T$ |
| | 0.8 | ⟨129.5760, 16.44⟩ | $(0\ 0\ 1\ 0\ 0\ 0\ 1\ 0\ 0\ 1\ 0\ 1\ 1\ 0\ 1\ 0\ 1\ 1)^T$ |

Crisp equivalent of b-$(\alpha, \beta)$ RFQMSTP, i.e., problems (24) and (25) are solved at a specific confidence levels of $\alpha$ and $\beta$. However while taking decisions a decision maker may alter the predetermined confidence levels under different dicision making situations. Such alterations in decisions have an impact on the solutions of the problem. In order to observe how the compromised solution changes for different confidence levels, we have perform a sensetivity analysis of the compromised solution of $G$, shown earlier in, Fig. 5, at different pairs of credibility ($\beta$) and trust level ($\alpha$) which are shown in Table 4. For this purpose, the compromised decision vector for a unique combination of credibility and trust measure are obtaind by solving problems (24) and (25) using epsilon-constrained method.

We have also simulated the models of problem (24) and problem (25) for larger graph instances. The nomenclature of all the random instances is of the form QMST_n_m, where n is the number of vertices and m is the number of edges for the graph $G$ and there exist a quadratic weight between any pair of edges in G. Both the linear and quadratic weights coefficients of the problems are generated randomly. The models of such problems are solved with MOGAs: NSGA-II [27] and MOCHC [23] and the corresponding performance metrics, e.g., hypervolume (HV) [13], spread (Sp) [4], generational distance (GD) [11], inverted generational distance (IGD) [11] and epsilon (E) [14] are studied to make a comparison between the MOGAs. In this context, it is quite imperative to note that there does not exist any set of optimal or nondominated solutions, i.e., Pareto front of any rough fuzzy bi-objective quadratic minimum spanning problem in the literature. Henceforth, the performance metrics for the MOGAs are measured based on a reference front which is constructed by gathering the best nondominated solutions from every independent execution of these MOGAs. The parameter settings for NSGA-II [23] and MOCHC [27] used to optimize the b-$(\alpha, \beta)$RFQMSTP model for the corresponding 5 rough fuzzy instances, generated randomly are as follows.

- For NSGA-II [23]: Population size = 100, function evaluations = 50,000, crossover probability = 0.9, mutation probability = 0.03.
- For MOCHC [32]: Population size = 100, function evaluations = 50,000, crossover probability = 0.9, incest threshold value = 0.25, persevered population = 0.05.

Both the MOGAs, i.e., NSGA-II and MOCHC are executed for 100 times for each of the rough fuzzy instances. For each execution, some performance metrics, e.g., HV, Sp, GD, IGD and E are evaluated with respect to the optimized solutions obtained after 50,000 function evaluations and the generated reference front. Due to the stochastic nature of the MOGAs, the measures of statistical dispersion among the performance metrics are studied which are listed in Table 5 through Table 8. Table 5 and Table 6 contain the statistical measures of the performance metrics when the credibility level ($\beta$) and the trust level ($\alpha$) for all the instances are set to $\alpha = 0.9$ and $\beta = 0.8$ respectively. Table 7 and Table 8 show the statistical measures of the performance metrics when $\alpha = 0.9$ and $\beta = 0.8$ respectively.

Considering the data in Table 4 through Table 7, we observe that out of 5 rough fuzzy instances, generated randomly, MOCHC [23] proves to be better than NSGA-II [27] for 3 instances. We also notice that for all the cases MOCHC emerges as better in terms of Sp [3] for every instances which means the diversity preservation among the solutions is more efficient in MOCHC [23] compared to NSGA-II [27]. For the credibility level 0.4 we noticed that MOCHC [23] is better in terms of GD [11] and IGD [11] for most of the instances whereas NSGA-II [26] is better in terms of HV [13], and E [14] for major number of instances. For the credibility level 0.8 we observe that HV improves for most of the instances for MOCHC [23] with respect to NSGA-II [27] but GD [11], IGD [11] and E [14] becomes more efficient for NSGA-II [27] for a large no of instances. However, considering all the performance matrices, MOCHC [23] outperforms NSGA-II [26] for 3 out of 5 random instances for two different combinations of at credibility ($\beta$) and trust ($\alpha$) levels, i.e., for $\beta = 0.4, \alpha = 0.9$ and $\beta = 0.8, \alpha = 0.9$ respectively.

The Java based metaheuristic framework jMETAL 4.5 [22] has been used for simulation of the MOGAs. NSGA-II [27] and MOCHC [27] have been executed for 50,000 function evaluations. Both the pairs of $\beta = 0.4, \alpha = 0.9$ and $\beta = 0.8, \alpha = 0.9$ are considered for problems (24) and (25) respectively which makes 2 unique combinations for the experimental analysis. One combination was set with a credibility level and trust level of 0.4 and 0.9 respectively. The other one was set with a credibility and trust levels of 0.8 and 0.9 correspondingly. Each simulation configuration was repeated 100 independent times, for the stochastic behaviour of the MOGAs.

Table 5: Mean and standard deviation of HV, Sp, GD, IGD and E after 100 runs of NSGA-II and MOCHC for credibility level 0.4 and trust level 0.9

| MOGAs | Rough Fuzzy Instances generated randomly | HV | | Sp | | GD | | IGD | | E | |
|---|---|---|---|---|---|---|---|---|---|---|---|
| | | Mean | S.D. | Mean | S.D. | Mean | S.D. | Mean | S.D. | Mean | S.D. |
| | | $\beta = 0.4, \alpha = 0.9$ | | | | | | | | | |
| *NSGA-II* | QMST_10_30 | **7.66E-1** | 3.6E-2 | 1.58E+0 | 4.7E-2 | 3.13E-3 | 1.6E-3 | 8.63E-4 | 3.1E-4 | 2.59E+1 | 1.4E+1 |
| | QMST_20_70 | **7.63E-1** | 3.3E-2 | 1.43E+0 | 4.4E-2 | 4.04E-3 | 1.4E-3 | **8.49E-4** | 3.2E-4 | **3.58E+1** | 1.2E+1 |
| | QMST_30_120 | **7.42E-1** | 2.1E-2 | 1.10E+0 | 5.8E-2 | 3.02E-3 | 1.3E-3 | **5.20E-4** | 1.7E-4 | **1.06E+2** | 4.4E+1 |
| | QMST_40_170 | 7.91E-1 | 2.0-E-2 | 1.16E+0 | 5.4E-2 | 3.39E-3 | 1.7E-3 | 2.48E-3 | 5.8E-4 | **1.02E+2** | 4.1E+1 |
| | QMST_50_220 | 8.01E-1 | 3.9E-2 | 1.18E+0 | 5.7E-2 | **1.21E-3** | 0.9E-3 | 1.45E-2 | 3.1E-3 | **1.65E+2** | 5.6E+1 |
| *MOCHC* | QMST_10_30 | 7.38E-1 | 3.5E-2 | **6.45E-1** | 4.5E-2 | **1.20E-3** | 0.9E-3 | **5.76E-4** | 3.5E-4 | **1.77E+1** | 1.7E+1 |
| | QMST_20_70 | 7.15E-1 | 3.6E-2 | **6.47E-1** | 5.8E-2 | **2.97E-3** | 1.5E-3 | 1.17E-3 | 4.8E-4 | 4.84E+1 | 2.0E+1 |
| | QMST_30_120 | 7.19E-1 | 2.5E-2 | **6.92E-1** | 4.6E-2 | **2.97E-3** | 1.8E-3 | 6.46E-4 | 2.3E-4 | 1.37E+2 | 4.0E+1 |
| | QMST_40_170 | **7.98E-1** | 1.9-E-2 | **7.35E-1** | 4.9E-2 | **3.81E-3** | 2.0E-3 | **2.42E-3** | 4.7E-4 | 1.18E+2 | 4.0E+1 |
| | QMST_50_220 | **8.11E-1** | 2.4E-2 | **7.26E-1** | 4.1E-2 | 1.25E-3 | 0.3E-3 | **1.13E-2** | 2.3E-3 | **1.92E+2** | 5.8E+1 |

Table 6: Median and I.Q.R. of HV, Sp, GD, IGD and E after 100 runs of NSGA-II and MOCHC for credibility level 0.4 and trust level 0.9

| MOGAs | Rough Fuzzy Instances generated randomly | HV | | Sp | | GD | | IGD | | E | |
|---|---|---|---|---|---|---|---|---|---|---|---|
| | | Median | I.Q.R. | Median | I.Q.R. | Median | I.Q.R. | Median | I.Q.R. | Median | I.Q.R. |
| | | $\beta = 0.4, \alpha = 0.9$ | | | | | | | | | |
| *NSGA-II* | QMST_10_30 | **7.91E-1** | 4.7E-2 | 1.57E+0 | 4.7E-2 | **1.34E-3** | 5.8E-3 | 8.20E-4 | 4.3E-4 | 3.26E+1 | 2.9E+1 |
| | QMST_20_70 | **7.76E-1** | 3.0E-2 | 1.42E+0 | 8.7E-2 | **3.21E-3** | 3.3E-3 | **8.47E-4** | 4.2E-4 | **3.68E+1** | 2.1E-1 |
| | QMST_30_120 | **7.47E-1** | 2.6E-2 | 1.11E+0 | 1.0E-2 | 2.65E-3 | 1.4E-3 | **4.98E-4** | 1.7E-4 | **1.04E+1** | 2.1E+1 |
| | QMST_40_170 | 7.93E-1 | 2.3-E-2 | 1.17E+0 | 3.4E-2 | 5.32E-3 | 1.9E-3 | 2.40E-3 | 6.0E-4 | 9.99E+1 | 6.9E+1 |
| | QMST_50_220 | 8.09E-1 | 3.6E-2 | 1.19E+0 | 1.4E-2 | 1.17E-2 | 5.1E-3 | 1.38E-2 | 3.3E-3 | **1.69E+2** | 7.3E+1 |
| *MOCHC* | QMST_10_30 | 7.28E-1 | 3.5E-2 | **6.45E-1** | 5.2E-2 | 1.20E-3 | 1.1E-3 | **7.68E-4** | 6.0E-4 | **5.50E+0** | 3.3E+1 |
| | QMST_20_70 | 7.24E-1 | 4.6E-2 | **6.55E-1** | 1.1E-2 | 1.56E-2 | 7.9E-3 | 1.08E-3 | 5.2E-4 | 4.71E+1 | 3.0E+1 |
| | QMST_30_120 | 7.17E-1 | 2.5E-2 | **6.96E-1** | 3.9E-2 | **2.61E-3** | 1.7E-3 | 6.02E-4 | 2.9E-4 | 1.43E+2 | 5.7E+1 |
| | QMST_40_170 | **7.97E-1** | 2.8-E-2 | **7.37E-1** | 5.1E-2 | **3.57E-3** | 2.5E-3 | **2.39E-3** | 8.1E-4 | 1.15E+2 | 6.1E+1 |
| | QMST_50_220 | **8.14E-1** | 3.5E-2 | **7.49E-1** | 6.3E-2 | **1.14E-2** | 8.0E-3 | **1.11E-2** | 3.1E-3 | **1.90E+2** | 7.3E+1 |

Table 7: Mean and standard deviation of HV, Sp, GD, IGD and E after 100 runs of NSGA-II and MOCHC for credibility level 0.8 and trust level 0.9

| MOGAs | Rough Fuzzy Instances generated randomly | $\beta = 0.8, \alpha = 0.9$ | | | | | | | | | |
|---|---|---|---|---|---|---|---|---|---|---|---|
| | | HV | | Sp | | GD | | IGD | | E | |
| | | Mean | S.D. | Mean | S.D. | Mean | S.D. | Mean | S.D. | Mean | S.D. |
| NSGA-II | QMST_10_30 | **7.92E-1** | 2.0E-2 | 1.39E+0 | 4.1E-2 | **1.34E-3** | 1.8E-3 | **3.34E-4** | 1.1E-4 | **2.55E+1** | 1.1E+1 |
| | QMST_20_70 | 7.72E-1 | 3.2E-2 | 1.26E+0 | 4.8E-2 | **2.18E-3** | 2.2E-3 | 5.42E-4 | 2.2E-4 | **4.57E+1** | 2.4E+1 |
| | QMST_30_120 | **7.93E-1** | 1.7E-2 | 9.90E-1 | 6.9E-2 | 2.50E-3 | 9.6E-4 | **7.95E-4** | 1.8E-4 | **1.25E+2** | 4.7E+1 |
| | QMST_40_170 | 7.41E-1 | 2.5-E-2 | 1.10E+0 | 7.7E-2 | 3.16E-3 | 1.5E-3 | **7.55E-4** | 2.4E-4 | **1.11E-+2** | 5.0E+1 |
| | QMST_50_220 | 7.92E-1 | 2.9E-2 | 1.17E+0 | 8.7E-2 | **6.59E-3** | 2.4E-3 | **7.07E-3** | 1.8E-3 | 1.80E+2 | 6.6E+1 |
| MOCHC | QMST_10_30 | 7.82E-1 | 1.8E-2 | **6.56E-1** | 3.5E-2 | 3.25E-3 | 3.0E-3 | 3.82E-4 | 1.1E-4 | 3.10E+1 | 8.5E+0 |
| | QMST_20_70 | **7.88E-1** | 2.6E-2 | **6.75E-1** | 5.0E-2 | 5.05E-3 | 3.7E-3 | **3.82E-4** | 1.7E-4 | 6.04E+1 | 2.0E+1 |
| | QMST_30_120 | 7.83E-1 | 1.7E-2 | **7.46E-1** | 4.2E-2 | **1.65E-3** | 6.4E-4 | 1.15E-3 | 3.4E-4 | 1.56E+2 | 3.7E+1 |
| | QMST_40_170 | **7.47E-1** | 3.3-E-2 | **6.89E-1** | 4.8E-2 | **2.87E-3** | 2.4E-3 | 7.62E-4 | 3.5E-4 | 1.30E+2 | 5.4E+1 |
| | QMST_50_220 | **7.99E-1** | 3.9E-2 | **7.28E-1** | 4.8E-2 | 8.13E-3 | 4.1E-3 | 8.55E-3 | 2.7E-3 | **1.73E+2** | 6.3E+1 |

Table 8: Median and I.Q.R of different performance metrics after 100 runs of NSGA-II and MOCHC for credibility level 0.8 and trust level 0.9

| MOGAs | Rough Fuzzy Instances generated randomly | $\beta = 0.8, \alpha = 0.9$ | | | | | | | | | |
|---|---|---|---|---|---|---|---|---|---|---|---|
| | | HV | | Sp | | GD | | IGD | | E | |
| | | Median | I.Q.R. | Median | I.Q.R. | Median | I.Q.R. | Median | I.Q.R. | Median | I.Q.R. |
| NSGA-II | QMST_10_30 | 7.80E-1 | 1.4E-2 | 1.38E+0 | 5.0E-2 | **9.34E-4** | 6.4E-4 | **3.26E-4** | 1.0E-4 | 3.50E+1 | 1.2E+1 |
| | QMST_20_70 | 7.74E-1 | 2.8E-2 | 1.26E+0 | 6.2E-2 | **1.41E-3** | 1.3E-3 | **3.47E-4** | 1.9E-4 | **4.52E+1** | 3.0E-1 |
| | QMST_30_120 | **7.95E-1** | 2.3E-2 | 9.97E-1 | 1.0E-1 | 2.35E-3 | 1.2E-3 | **8.02E-4** | 2.9E-4 | **1.26E+2** | 7.3E+1 |
| | QMST_40_170 | 7.43E-1 | 4.5-E-2 | 1.10E+0 | 1.1E-1 | 2.90E-3 | 2.0E-3 | 7.17E-4 | 3.2E-4 | **1.04E+2** | 7.6E+1 |
| | QMST_50_220 | 7.94E-1 | 3.9E-2 | 1.17E+0 | 1.2E-1 | **6.27E-3** | 3.5E-3 | **6.78E-3** | 2.4E-3 | 1.80E+2 | 9.3E+1 |
| MOCHC | QMST_10_30 | **7.91E-1** | 2.4E-2 | **6.66E-1** | 4.7E-2 | 2.24E-3 | 2.6E-3 | 3.74E-4 | 1.2E-4 | **3.11E+1** | 2.2E+1 |
| | QMST_20_70 | **7.87E-1** | 4.4E-2 | **6.75E-1** | 7.3E-2 | 3.80E-3 | 5.0E-3 | 5.05E-4 | 2.6E-4 | 5.84E+1 | 2.7E+1 |
| | QMST_30_120 | 7.81E-1 | 2.6E-2 | **7.44E-1** | 5.2E-2 | **1.48E-3** | 7.7E-4 | 1.09E-3 | 3.4E-4 | 1.57E+2 | 4.7E+1 |
| | QMST_40_170 | **7.44E-1** | 3.5-E-2 | **6.96E-1** | 5.7E-2 | **2.04E-3** | 2.2E-3 | **6.91E-4** | 4.1E-4 | 1.25E+2 | 7.8E+1 |
| | QMST_50_220 | **7.97E-1** | 3.7E-2 | **7.37E-1** | 6.5E-2 | 8.02E-3 | 6.2E-3 | 8.29E-3 | 3.5E-3 | **1.69E+2** | 8.4E+1 |

## 8. Conclusion

In this article we have presented the rough fuzzy quadratic minimal spanning tree problem which is modelled as a bi-objective problem using rough fuzzy chance-constrained programming technique.

Till now no network model has been extended under such an uncertain environment in the literature where the weight parameters of the network are essentially rough fuzzy in nature. The proposed model is eventually solved by epsilon-constraint method, NSGA-II and MOCHC. This model is capable to capture hybrid uncertainty expressed in the form of rough fuzzy variables.

The rough fuzzy chance constrained model for quadratic spanning tree is also solved for 5 randomly generated rough fuzzy graphical instances using two MOGAs: NSGA-II and MOCHC. A comparative analysis between both the MOGAs is also performed in terms of different performance metrics.

The rough fuzzy chance constrained programming model can also be extended to other network problems like shortest path, network flow, transportation model, in near future.